\documentstyle[11pt,amssymb]{article}
\input{epsf}

\title{ Exponential equations for  the quantum ``$az+b$'' group.}

\author{Ma{\l }gorzata Rowicka - Kudlicka\thanks{Supported by KBN grant 
No 5 PO3A 036 18}\\
Institute of Mathematics, Polish Academy of Sciences,\\
\'Sniadeckich 8, 00-950 Warszawa, Poland\\
e-mail:$\;$rowicka@fuw.edu.pl}

\begin{document}

\maketitle

\newcommand{\mqed}{\nopagebreak\centerline{\hfill
\raisebox{3.5ex}[0ex][0ex]{$\Box$}}\\}
\newcommand{\faz}{{\rm Phase\  }}
\newcommand{\ci}{continuous}
\newcommand{\ru}{ unitary representation}
\newcommand{\Qa}{ Quantum 'ax+b' group}
\newcommand{\qa}{ quantum 'ax+b' group}
\newcommand{\MU}{ multiplicative unitary}
\newcommand{\slw}{S.L. Woronowicz}
\newcommand{\tr}{\Delta}
\newcommand{\cH}{{\cal H}}
\newcommand{\K}{{\cal K}}
\newcommand{\cK}{{\cal K}}
\newcommand{\sH}{_{\cal H}}
\newcommand{\sK}{_{\cal K}}
\newcommand{\h}[1]{\hat{#1}}
\newcommand{\lh}{L({\cal H})}
\newcommand{\rf}[1]{{\rm (\ref{#1})}}
\newcommand{\ch}{ {\cal C}(H)}
\newcommand{\chn}{ {\cal C}(H)^N}
\newcommand{\ov}{\overline}
\newcommand{\fh}{F_{\hbar}} 
\newcommand{\hb}{\hbar}
\newcommand{\vt}{V_{\theta}} 
\newcommand{\za}{-\!\!\circ} 
\newcommand{\R}{{\Bbb R}}
\newcommand{\C}{{\Bbb C}}
\newcommand{\Z}{{\Bbb Z}}
\newcommand{\ro}{\rho}
\newcommand{\si}{\sigma}
\newcommand{\be}{\beta}
\newcommand{\de}{\delta}
\newcommand{\ga}{\gamma}
\newcommand{\ta}{\tau}
\newcommand{\N}{{\Bbb N}}
\newcommand{\M}{{\rm M}}
\newcommand{\B}{{\rm B}}
\newcommand{\Lin}{{\rm L}}
\newcommand{\pod}{{\rm d}}
\newcommand{\iz}{\cong}
\newcommand{\eps}{\epsilon}
\newcommand{\fil}{\varphi}
\newcommand{\la}{\lambda}
\newcommand{\qed}{ $\Box$}
\newcommand{\Ci}{C_{\infty}}
\newcommand{\Cir}{C_{\infty}(\R)}
\newcommand{\eh}{e^{\frac{i\hbar}{2}} }
\newcommand{\ehm}{e^{-\frac{i\hbar}{2}} }
\newcommand{\Cgr}{C^{\infty}(\R)}
\newcommand{\Cg}{C^{\infty}}
\newcommand{\Cor}{C_{\rm o}(\R)}
\newcommand{\Co}{C_{\rm o}}
\newcommand{\Cog}{C_{\rm bounded}}
\newcommand{\Cogr}{C_{\rm b}(\R)}
\newcommand{\Lk}{L^{2}}
\newcommand{\Lkr}{L^{2}(\R)}
\newcommand{\po}{\hat{p}}
\newcommand{\qo}{\hat{q}}
\newcommand{\xo}{\hat{x}}
\newcommand{\csta}{$C^{*}$-}
\newcommand{\cstal}{$C^{*}$-algebra \ }
\newcommand{\te}{\otimes}
\newcommand{\ad}{{\rm ad}}
\newcommand{\id}{{\rm id}}
\newcommand{\Mor}{{\rm Mor}}
\newcommand{\Rep}{{\rm Rep}}
\newcommand{\spe}{{\rm Sp }}
\newcommand{\sign}{{\rm sign }\;}
\newcommand{\whe}{\hspace*{5mm}\mbox{\rm where}\hspace{5mm}}
\newcommand{\mand}{\hspace*{5mm} {\rm and} \hspace{5mm}}
\newcommand{\moraz}{\hspace*{5mm} {\rm and} \hspace{5mm}}
\newcommand{\af}{\hspace*{1mm} {\bf \eta} \hspace{1mm}}
\newcommand{\od}{\hspace*{5mm}}
\newcommand{\fu}{{\cal F}}
\newcommand{\fuod}{{\cal F}^{-1}}
\newcommand{\mlot}{\mbox{$\hspace{.5mm}\bigcirc\hspace{-3.mm}
\raisebox{-.7mm}{$\top$}\hspace{1mm}$}}
\newcommand{\Zak}{\mbox{$-\hspace{-2pt}\comp\,$}}
\newcommand{\dow}{{\bf Proof: }}
\newcommand{\ut}{\cong}
\newcommand{\Sp}{{\rm sp}}
\newcommand{\infi}{\infty}
\newcommand{\tend}{\rightarrow}
\newcommand{\impl}{\Rightarrow}
\newcommand{\limn}{\lim_{n\rightarrow\infty}}
\newcommand{\limk}{\lim_{k\rightarrow\infty}}
\newcommand{\limt}{\lim_{t\rightarrow\infty}}
\newcommand{\limx}{\lim_{|x|\rightarrow\infty}}
\newcommand{\lime}{\lim_{\eps\rightarrow 0}}
\newcommand{\Lj}{L^{1}}
\newcommand{\Lp}{L^{P}}
\newcommand{\Ls}{L^{S}}
\newcommand{\Li}{L^{\infty}}
\newcommand{\ljr}{{\cal L}^{1}(\R)}
\newcommand{\lk}{ l^{2}}
\newcommand{\lpr}{{\cal L}^{P}(\R)}
\newcommand{\lir}{{\cal L}^{\infty}(\R)}
\newcommand{\Ljr}{L^{1}(\R)}
\newcommand{\Lpr}{L^{P}(\R)}
\newcommand{\Lsr}{L^{S}(\R)}
\newcommand{\Lir}{L^{\infty}(\R)}
\newcommand{\lkn}{ l^2(\N)}
\newcommand{\lp}{{\cal L}^{P}}
\newcommand{\li}{{\cal L}^{\infty}}
\newcommand{\Ljx}{L^{1}(\X)}
\newcommand{\Lkx}{L^{2}(\X)}
\newcommand{\Lpx}{L^{p}(\X)}
\newcommand{\Lqx}{L^{q}(\X)}
\newcommand{\Lix}{L^{\infty}(\X)}
\newcommand{\ljx}{{\cal L}^{1}(\X)}
\newcommand{\lkx}{{\cal L}^{2}(\X)}
\newcommand{\lpx}{{\cal L}^{p}(\X)}
\newcommand{\lix}{{\cal L}^{\infty}(\X)}
\newcommand{\ilkC}{C_{\infty}(\R)\otimes_{C}C_{\infty}(\R)}
\newcommand{\ilk}{\otimes_{\cal C}}
\newcommand{\zw}{CB(\Lkr)}
\newcommand{\zwg}{CB\left(\,L^{2}(G)\,\right)}
\newcommand{\ogr}{B(L^2(\R))}
\newcommand{\hs}{HS(H)}
\newcommand{\spl}{\star}
\newcommand{\milk}{$(i_{1},i_{2},\zw)$}
\newcommand{\mczw}{{\rm Mor(\Cir,\zw)}}
\newcommand{\res}{{\rm Res}}
\newcommand{\pra}{\mbox{${\rm Proj_1}$}}
\newcommand{\prb}{\mbox{${\rm Proj_2}$}}
\newcommand{\Morc}{{\rm Mor}_{\katC}}
\newcommand{\recogr}{{\rm Rep(\Cir,\ogr)}}
\newcommand{\mi}{\hspace*{3mm} {\rm and} \hspace{3mm}}
\newcommand{\mor}{\hspace*{5mm} {\rm or} \hspace{5mm}}
\newcommand{\dla}{\hspace*{5mm} {\rm for} \hspace{5mm}}
\newcommand{\ja}{j_{1}}
\newcommand{\ib}{i_{2}}
\newcommand{\ia}{i_{1}}
\newcommand{\jb}{j_{2}}
\newcommand{\jc}{j_{3}}
\newcommand{\sz}{{\cal S}(\R)}
\newcommand{\gw}{$^*$}
\newcommand{\heps}{h_{\eps}} 
\newcommand{\czn}{\Co}
\newcommand{\bh}{B(H) }
\newcommand{\dlad}{for any }
\newcommand{\dlak}{for each }
\newcommand{\Gn}{{\cal G}_n}
\newcommand{\G }{{\rm G}_n}
\newcommand{\rt}{{\bf T}^1}
\newcommand{\jd}{{\bf T}^2}
\newcommand{\Tau}{{\cal T}}


\newcommand{\bfa}{\begin{fakt}}\newcommand{\efa}{\end{fakt}}
\newcommand{\ble}{\begin{lem}}\newcommand{\ele}{\end{lem}}
\newcommand{\bst}{\begin{stw}}\newcommand{\est}{\end{stw}}
\newcommand{\bde}{\begin{defi}}\newcommand{\ede}{\end{defi}}
\newcommand{\bwn}{\begin{wn}}\newcommand{\ewn}{\end{wn}}
\newcommand{\buw}{\begin{uwaga}}\newcommand{\euw}{\end{uwaga}}
\newcommand{\bdy}{\begin{dygresja}}\newcommand{\edy}{\end{dygresja}}
\newcommand{\bwa}{\begin{warning}}\newcommand{\ewa}{\end{warning}}
\newcommand{\bpr}{\begin{przy}}\newcommand{\epr}{\end{przy}}
\newcommand{\btw}{\begin{tw}}\newcommand{\etw}{\end{tw}}
\newcommand{\beq}{\begin{equation}}\newcommand{\eeq}{\end{equation}}
\newcommand{\bit}{\begin{itemize}}\newcommand{\eit}{\end{itemize}}
\newcommand{\bq}{\begin{quote}}\newcommand{\eq}{\end{quote}}
\newcommand{\ba}{\begin{array}}\newcommand{\ea}{\end{array}}


\newtheorem{defi}{Definition}[section]
\newtheorem{wn}[defi]{Observation}
\newtheorem{tw}[defi]{Theorem}
\newtheorem{lem}[defi]{Lemma}
\newtheorem{fakt}[defi]{Corollary}
\newtheorem{stw}[defi]{Proposition}
\newtheorem{przy}[defi]{Example}
\newtheorem{uwaga}[defi]{Remark}
\newtheorem{warning}[defi]{Warning}


\begin{abstract} 
We consider quantum group theory on the Hilbert space level. 
We find all solutions for scalar and general exponential equations 
for the quantum ``$az+b$'' group. It turns out that there is a 
simple formula for all of them involving the quantum 
exponential function $F_N$.
The very interesting theorem we prove by the way  
 is the one on the existence of normal extension 
of certain sum of normal operators. 
          
To put it differently,   we find all unitary 
representations 
of the braided quantum group related 
to the quantum ``$az+b$'' group. 
 This is the most difficult 
result needed to classify all unitary representations
 of the quantum ``$az+b$'' group. 
Eventually this  enables us to give a formula 
for all unitary representations of the quantum ``$ax+b$'' 
group in our next paper \cite{paper4}.

{\bf key words:} unbounded operators -- Hilbert space

{\bf MSC-class:} 20G42 (Primary), 47B25 (Secondary).
\end{abstract}

\section{Introduction} 

To explain  what is going on 
 in this paper  let us use an analogy with the classical case.
One of the goals of the classical group theory 
 is to find all unitary representations of the group considered. 
For example, by SNAG theorem, we know that $U$ is a  
(strongly continuous) unitary representation
 of the group  $\R^2$ acting on Hilbert space $\cH$ iff  $\Ci(\R^2)$
there exists  a pair of strongly 
commuting selfadjoint operators $(a,b)$ acting on $\cH$ 
such  
 that for any $(x,y)\in \R^2$  we have
\[U(x,y)=e^{ixa+iyb}\ .\]
It means that all unitary representations of $\R^2$ 
 are ``numbered'' by elements from the set of all 
pairs of strongly commuting selfadjoint operators.

On the other hand, every such pair gives rise, through the functional 
calculus of normal operators, to a representation of the algebra of all 
continuous vanishing at infinity functions on the group dual to 
$\R^2$, i.e. in this case to representations of 
 $\Ci(\hat{\R^2})=\Ci(\R^2)$. This phenomenon, 
i.e. correspondence between unitary representations of the locally compact 
group $G$ and representations of the algebra  $\Ci(G)$, is known as 
the Pontryagin duality.   
On the other hand, the exponential function is a solution 
 of the exponential equation
\beq
\label{rownwykl}
F(x+y)=F(x)F(y),
\eeq
 where $x,y\in\R^2$. Moreover, if we look for a solution 
such that $F$ is measurable and $|F(x)|$ has a modulus 1 for every $x$,
  all the  solutions of \rf{rownwykl}
 is given by  the formula
\[F(x,y)=e^{ixa+iyb}\ ,\]
where $(a,b)\in\hat{\R^2}$.

The solution of the general exponential equation, i.e. with 
 the one unitary-operator-valued function $F$, is a direct integral 
of the solutions of the scalar case. So once one knows solutions 
of the general exponential equation, one knows also all unitary 
representations of the group involved, in this case $\R^2$.

One can also consider equation \rf{rownwykl} in more general setting.
One can allow $x$ and $y$ to be ``coordinates'' on the two copies 
of a space $G$, classical or quantum. They are usually sets of operators 
 acting on Hilbert space and satisfying certain conditions. 
8888

This paper is very similar in spirit to the our forthcoming paper 
\cite{paper1}, where braided quantum groups related to the quantum 
``ax+b'' group were studied and their unitary representations classified.

In this paper we consider the quantum group theory  on the 
Hilbert space level. The quantum ``az+b'' group was constructed
 recently  by 
S.L. Woronowicz in \cite{az+b}. 
It is the natural deformation of the group of affine transformations of the 
 complex plane, with the deformation parameter $q$ 
 being an even primitive root of unity (for details see cite{az+b}). 
However, what we are mainly interested in in this paper is a braided 
quantum group related to the quantum  ``az+b'' group.

Let us begin with an explanation, what we mean by a quantum group there.

In fact, the definition of locally compact quantum group, 
and such is ``az+b'',
is still under construction  \cite{vaeskrotki},  however  one knows 
approximately what a quantum group should be.

We like most the approach using operator domain and operator 
function (described in \cite{wunb,oper,paper1}). 
We believe it is a  very insightful  one.
What we say below is not really necessary to understand the paper,
 but we hope it will be useful to understand idea that lie behind and 
the connections with the quantum ``az+b'' group.

We only outline general ideas here, for more detailed treatment we 
refer the Reader to \cite{wunb,oper,paper1}.

The description of an operator domain is similar to the global one of 
a manifold, where coordinates and relation satisfied by them are given.
In our noncommutative case, the coordinates are closed (so in 
general unbounded) operators and the relations are arbitrary commutation 
rules  that are invariant with respect to unitary transformations and 
direct sum decomposition. For this invariance S.L. Woronowicz, 
to whom this idea is due, coined 
 a term ``respecting symmetry of the Hilbert space''.

Easy example of an operator domain is the one crucial in this paper:
 an operator domain  $D$ related to the quantum ``az+b'' group.
The relations in this case are
\[RR^*=R^*R \mi\spe R\subset \ov\Gamma \;\,\]
where $\Gamma$ is the multiplicative subgroup of $\C\setminus\{0\}$ 
given by \rf{gamma}.
The unbounded operator $R$ entering descriptions of the operator domain 
$D$ 
 can be thus thought of as a ``coordinate on a  quantum space''. 
Observe, that this space is entirely classical, it can be identified 
with $\ov{\Gamma}$.

The operator functions can be thought of as a  recipe 
what to do with a $N$-tuple of closed 
operators $(a_1, a_2,..., a_N)$ to obtain another closed operator 
 $F(a_1, a_2,..., a_N)$.
An operator map are similar to operator functions, the only difference is 
they may transform $N$-tuple of closed 
operators $(a_1, a_2,..., a_N)$ to obtain another $k$-tuple of closed 
operators  $(F_1(a_1, a_2,..., a_N),F_1(a_1, a_2,..., a_N),)$

 Let $G$ be an operator domain and let $G\times G$ denote an operator 
domain
\[G\times G:=\{ (x,y)\;|\; x,y\in G \mi xy=yx\}\]
Let $\cdot$ be an operator map
$$\cdot{ } \;\;\;:G\times G\ni (x,y)\rightarrow xy\in G$$.

Loosely speaking, a quantum  group $G$  is such an operator 
domain $G$ equipped with an associative operator map $\cdot$.

\bpr[Quantum  ``$a{z}+b$'' at roots of unity] 
$$G=\left\{\;(a,b) \;
:\;
\ba{c} \; \\
aa^*=a*a,\; bb^*=b^*b,\\
 a \;\;\mbox{\rm invertible} \\
a b=qb a\\
\spe a,\spe b\subset \ov\Gamma ,
\ea
\right\}
 $$
where $q$ in an even primitive root of unity and $\Gamma$ 
is given by \rf{gamma}.

Group operation in  $G $ is given by
\[
\cdot:
G\times 
G\ni ((a_1,b_1),
(a_2,b_2))
\tend (a,b)\in G_{\rm H}\]
where 
\[a =a_1\te a_2\mi 
b=a_1\te  b_2{\dot{+}} b_1\te I,\]
\epr
where $\dot{+}$ denotes closure.
The spectral condition is there to make 
 sure that $G$ is closed under the operation $\cdot$.

 The  main difference between a braided quantum group and 
 a quantum group is that a group operation on a braided 
 quantum group $G$ is defined on  a smaller  operator domain 
\[G^2:=\{ (x,y)\;|\; x,y\in G \mi x,y \
\mbox{\rm satisfy certain relations}\}\;.\]
Usually we do not assume that operators from both copies 
 of $G$ commute, so in general a braided quantum group is 
not a quantum group.
A group operation on a braided quantum group $G$ should be the  
operator map 
\[\mlot:G^2\ni (x,y)\rightarrow x\mlot y\in G\;\] 
 which is associative.

\bpr[Braided quantum group $D$]
Let us define operator domains
\[D=\{R\;|\;RR^*=R^*R \mi\spe R\subset \ov\Gamma \;\;\}\]
and
\[\!\!D^2=
\{\left(R,S)\right)|\;
R,S\in D\;\;(\faz S)R=q R(\faz S)\;\]
\[\mbox{\rm and on  }(\ker S)^{\perp} \;
\mbox{\rm we have }\;
 |S|^{it}R=e^{-\frac{2\pi}{N}t}R|S|^{it} \;{\rm for any}
\; t\in\R\; \}\]
and $N=6,8,10,...$.

We define an operation  
$$\mlot_D\;:\; D^2\rightarrow D$$ 
by
\[R\mlot_D S
= R{\dot{+}}S.\]
Thus defined operation $\mlot_D$ is associative and 
$D$ with this operation  forms a  braided quantum  
group.

Moreover, let us observe that the braided quantum group $D$ 
is related to the quantum ``az+b'' group, denoted here by $G$, 
in the following way 
\[  (a,b)\in G
\Longleftrightarrow   
\left( \ba{c} (b,a)\in D^2\\
\ker a=\{0\} \ea\right)\ . \]
\epr

In this paper we find 
 all solutions of the general exponential equation for the quantum 
``az+b'' group.
In fact, we find  all the unitary representations 
of the braided quantum group $D$.
The main result is Theorem \ref{hip2}, which 
gives formula for all such representations.
This result is essential for classification of  
   all 
unitary representations 
  of the quantum  group "az+b", which is 
 achieved   in our forthcoming paper 
 \cite{paper4}.

The second important result in this paper is Proposition \ref{nieu}, 
which solves the problem of the existence of a normal extension 
of a sum $\mu R+S$, where $(R,S)\in D^2$ and $\mu \in \Gamma$.
We hope that this result will be useful in the construction of the quantum 
$GL(2,\C)$ group \cite{pgl2c}.

 In the remaining part of this section  we introduce   some non-standard
 notation and  notions used in this paper.

In Section \ref{sec:com} we introduce commutation rules related to the 
quantum ``az+b'' group. In the next section we discuss properties of 
pair of operators $(R,S)$ satisfying these commutation rules. In Section 
 \ref{sec:theq} we repeat the definition of the quantum exponential 
 function for the ``az+b'' group after \cite{az+b}. Then we investigate 
 in Section \ref{sec:nor} 
 the existence of normal extensions of 
$\mu R+S$, where $(R,S)\in D^2$ and $\mu \in \Gamma$.
In Section \ref{sec:sols} we give all solutions of the scalar exponential 
equation for the quantum ``az+b'' group and finally in Section 
\ref{sect:az2} we do the same for the general exponential equation.

In Appendix \ref{dodsuma} we prove the formula we use in Section 
\ref{sect:az2}.

\subsection{Notation}
We  denote Hilbert spaces by   
$ \cH$ and $\cK$, the set of all closed operators acting on $\cH$ by 
$\ch$, the set of bounded operators by  
$\B (\cH)$ and the sets of compact and unitary ones 
by $ CB (\cH)$
and ${\rm Unit}(\cH)$, respectively.
 The set of all continuous vanishing at infinity functions on a space $X$ will 
be denoted by $\Ci (X)$.
We consider only separable Hilbert spaces, usually infinite dimensional.
We  denote scalar product by
$(\cdot|\cdot)$ and it is antilinear in the first variable.
We  consider mainly unbounded linear operators.
All operators considered are densely defined. 
We  use functional calculus of normal operators 
 \cite{reedI,reedtf,  rudaf}.
We also use the symbol $\faz T$  for partial isometry  
  obtained from  polar decomposition of a normal  operator $T$.

We  use a non-standard, but very useful notation for orthogonal projections 
and their images  \cite{qexp}, as explained
 below. 
Let $a$ and $b$ be  strongly commuting selfadjoint  
operators  acting on  a Hilbert space  $\cH$. Then by spectral  theorem
there exists a common spectral measure $dE(\lambda)$ such that
\[
a={\int_{\R^2}} \lambda\, dE(\lambda,\mu),
\hspace{1cm}
b={\int_{\R^2}} \mu\, dE(\lambda,\mu)
.\]
For every complex measurable   function    
$f$ of two variables
\[
f(a,b)={\int_{\R^2}} f(\lambda,\lambda')\, dE(\lambda,\lambda').
\]

Let $f$  be a  logical sentence and let   $\chi(f)$
 be 0 if is false, and 1 otherwise.
If ${\cal R}$ is a binary relation on $\R$  
then $f(\lambda,\lambda')=\chi({\cal R}(\lambda,\lambda'))$ is   
 a  characteristic function of a set
$$\Delta=
\{(\lambda,\lambda')\in\R^2:{\cal R}(\lambda,\lambda')\}$$
 and assuming that  $\Delta$ is  measurable $f(a,b)=E(\Delta)$. 
From now on we will write $\chi({\cal R}(a,b))$
 instead of  $f(a,b)$:
\[
\chi({\cal R}(a,b))=
{\int_{\R^2}} \chi({\cal R}(\lambda,\lambda'))\, dE(\lambda,\lambda')
=E(\Delta).
\]
Image of this projector will be denoted  by  ${\cal H}(
{\cal R}(a,b))$, where `$\cH$' is a  Hilbert space, 
 on which operators  $a,b$ act.

Thus  we defined symbols  $\chi(a>b)$, 
 $\chi(a^2+b^2=1)$, \mbox{$\chi(a=1)$}, $\chi(b<0)$, $\chi(a\neq 0)$ 
 etc. They are  orthogonal projections  on appropriate 
  spectral subspaces. For example  ${\cH}(a=1)$ is  
 is an eigenspace of  operator  $a$ for eigenvalue 
 $1$ and   $\chi(a=1)$ is   orthogonal projector on   this  eigenspace.

Generally, whenever  $\Delta$ is a measurable subset of  $\R$,
then  ${\cH}(a\in\Delta)$ is   spectral subspace of an   operator 
$a$ corresponding to   $\Delta$ and  $\chi(a\in\Delta)$ is its   spectral projection.

Let $\qo$ and  $\po$   denote the position  and momentum operators in 
 Schr\"odinger representation, i.e. we set $\cH=\Lkr$. Then 
 the domain of $\qo$
$$D(\qo)=\{\psi\in \Lkr\;:\; \int_{\R}x^2|\psi(x)|^2dx<\infty\;\}$$
and $\qo$ is multiplication by  coordinate operator on that domain
$$(\qo \psi)(x)=x\psi(x).$$
The domain of  $\po$ consists of all distributions from   $\Lkr$ such that  
$$D(\po)=\{\psi\in \Lkr\;:\; \psi'\in\Lkr\;\} $$
and for any $\psi\in D(\po)$ 
$$(\po f)(x)=\frac{\hbar}{i}\frac{d f(x)}{dx}
,$$
where $-\pi<\hb<\pi$.

\section{Commutation rules related to the quantum "az+b" group}
\label{sec:com}
Let
\beq
\label{q}
q=e^{\frac{2\pi i}{N}}\ ,
\eeq
where $N$ is an even number and  $N\geq 6$, 
i.e. $q$ is a primitive root of unity: $q^N=1$.
Let us introduce notation
\beq
\label{hb2}
\hb=\frac{2\pi}{N}\ .
\eeq
Note that $\hb<\pi$ and  $q=e^{i\hb}$.

The assumption that  $N$ is an even number was made by  
S.L. Woronowicz in \cite{az+b}, where the quantum  
"az+b" was constructed. We need the assumption that  $N\geq 6$  
 to use formulas  (1.31), (1.32) and (1.34) from \cite{qexp}
(or equivalently (1.10), (1.11) and (1.12) from \cite{az+b}).
We will use these formula to derive relations 
we need later on.

Let
\beq
\label{gamma}
\Gamma =\bigcup_{k=0}^{N-1}q^k\;\R_+\ .
\eeq
It means that  $\Gamma$ is a multiplicative group. 

Let  $\ov\Gamma$ denote a closure of  $\Gamma$ in $\C$, i.e.
\beq
\label{gammaov}
\ov\Gamma =\Gamma \cup \{0\}\ .
\eeq
Our goal is to find an exponential equation
 for a pair of operators  $(R,S)$ acting on  Hilbert
 space $\cH$ and satisfying commutation relations 
 described in  \cite{az+b} and denoted by  
$(R,S)\in D_{\cH}$ there. 
We recall below definition of these relations.
\bde
\label{Dmu}
Let $\cH$ be a separable Hilbert space.
We say that closed  operator  $(R,S)\in D_{\cH}$ if
\bit
\item[1. ] $R,S$ are normal
\item[2. ] $\ker R=\ker S=\{0\}$ 
\item[3. ] $\spe R,\ \spe S\subset \ov\Gamma$
\item[4. ]\[({\rm Phase} S) \;R=q\;R \;({\rm Phase} S)\]
\[\mi |S|^{it}\;R=e^{-\frac{2\pi}{N}t}\;R\;|S|^{it}\ ,
\]
for any  $t\in\R$.
 \eit 
\ede
\buw
\label{azu1}
Condition  $4.$ can be written in equivalent form using  
 polar decomposition of operators   $R$ and  $S$
\[({\rm Phase} S) \;|R|=|R| \;({\rm Phase} S)\]
\[|S|\; ({\rm Phase} R)=({\rm Phase} R) \; |S|\]
\[({\rm Phase} S) ({\rm Phase} R)=q\;({\rm Phase} R) \;({\rm Phase} S) \]
and 
\beq
|S|^{it}\;|R|\;|S|^{-it}=e^{-\frac{2\pi}{N}t}\;|R| \ ,
\label{komut}
\eeq
for any $t\in\R$.

\euw

Observe also that if one performs analytical continuation and substitutes 
$t=i$ in \rf{komut}, one gets
$|S|\;|R|=q\;|R|\;|S|\ .$
Hence, using other equalities introduced in Remark
 \ref{azu1} we get
\[SR=q^2RS \mi SR^*=R^*S\ .\]
More precise derivation of these formulas can be found 
 in proof of 
Proposition 2.1  \cite{az+b}.

\section{Properties of  operators $(R,S)\in D_{\cH}$}
\label{sec:pro}
We prove now   analogues of  Proposition 3.1, Proposition 3.2 and 
Theorem 3.3  \cite{qexp}. 
The proofs below are modification of those presented in the cited
above paper.

Let us consider Hilbert space $\Lk (\Gamma,d\gamma)$, where
 $d\gamma$ is   Haar measure  
 of the group $\Gamma$, i.e.
\[\int_\Gamma f(\gamma)d\gamma =\sum_{k=1}^N\int_0^\infty f(q^kr)\frac{dr}
{r} ,\] 
 for any $f\in\Lk (\Gamma, d\gamma)$ and  $r\in\R_+$ and $k=0,1,...,N-1$. 
 
Let
\[\Lambda_k=\left\{z:\;\frac{2k\pi}{ N}<\arg z<\frac{2(k+1)\pi}{ N}\right\}\ \]
and
\beq
\label{dlamk}
\ov\Lambda_k=\left\{z:\;\frac{2k\pi}{ N}\leq\arg z\leq\frac{2(k+1)\pi}{ N}
\;\;\mbox{\rm or }\; z=0\ . \right\}\eeq
Then
\[\C=\Gamma\cup \bigcup_{k=0}^{N-1}\Lambda_k\ .\]
Let 
\[\Gamma_k=e^{\frac{2k\pi i}{ N}}\R_+\ .\]
Then
\[  \Gamma=\bigcup_{k=0}^{N-1}\Gamma_k\ .\]
\bde
\label{defHbound}
Let $H _{bounded}$ denote the set of all functions $f\in C(\Gamma)$, 
 such that there exists  a continuous and bounded  
function $\tilde{f}$ defined on  $\ov{\Lambda_o}\times q^{\Z}$ and such that 
\bit
\item[1.] for any  $k\in \Z$ function 
\[\Lambda_o\ni z\rightarrow \tilde{f}(z,q^k)\]
 is holomorphic
\item[2.] for any  $r\in\R_+$ and $k\in\Z$
\beq
\label{wl2}f(q^kr)=\tilde{f}(r,q^k)\ .
\eeq
\eit
\ede
If $f\in H _{bounded}$, it follows that there exist the described above 
  function $\tilde f$. We will use notation
\[f(q\cdot q^kr):=\tilde{f}(qr,q^k)\ .\]
\bst
\label{QE3.1}
Let $(R,S)\in D_{\cH}$ and $f\in H _{bounded}$. 
Then
\[f(R)D(S)\subset D(S)\]
and for any   $x\in D(S)$ we have  
\beq
\label{QE3.1wz}
f(q\cdot qR)Sx=Sf(R)x.
\eeq
\est
  \dow
Let $\Sigma$ be the  stripe   $\{\tau\in \C:0<\Im\tau<1\}$ and let  
$H(\Sigma)$ be the space of all  functions continuous on  
 $\ov\Sigma$
 and holomorphic in the interior of   $\Sigma$.
 Observe that space  $H (\Sigma)$  equipped with 
sup norm is a  Banach space.

For any  $\lambda\in\Gamma$ and $\tau\in\Sigma$ let us introduce 
 notation 
 \[\varphi_\lambda(\tau)=\tilde{f}\left(e^{\hb\tau}\lambda,q^k\right)\ ,\]
where $k=0,1,...,N-1$, and the relation between  $\tilde{f}$  and $f$ 
 is given by 
 \rf{wl2}.
 Then \mbox{$\varphi_\lambda\in H (\Sigma)$} and
$\|\varphi_\lambda\|\leq C$, where
 $C=\sup\{|f(r,q^k)|:r\in\Lambda_o \mi k\in \{0,1,...,N-1\}\}$. 
Therefore
\beq
\label{3.3}
\int_\Gamma \varphi_\lambda (\tau) d\mu(\lambda)\in H(\Sigma)
\eeq
for any (complex-valued) finite measure 
 $d\mu(\lambda)$ on $\Gamma$.

Let $dE_R(\lambda)$ be a spectral  measure of a normal operator  $R$,
 $x,y\in D(S)$ and 
$d\mu(\lambda)=(y|dE_R(\lambda)Sx)$. Then 
\[
\int_\Gamma \varphi_\lambda(\tau) d\mu(\lambda)=
(y|\tilde{f}(e^{\hbar\tau}|R|,\faz R)Sx).
\]
Moreover, by \rf{3.3} the map
\[
\Sigma\ni\tau\longrightarrow (y|\tilde{f}(e^{\hbar\tau}|R|,\faz R)Sx)\in\C
\]
 is a   function continuous on $\ov\Sigma$ and holomorphic inside 
 $\Sigma$.

Since  $x,y\in D(S)$, it follows that 
\[
\ov\Sigma\ni\tau\longrightarrow |S|^{1+i\tau}x\in \cH,
\]
\[
\ov\Sigma\ni\tau\longrightarrow |S|^{i\overline{\tau}}y\in \cH
\]
are continuous. The former function is    holomorphic, whereas  
 the latter is  antiholomorphic  on $\Sigma$. Therefore the  function  
\[
\ov\Sigma\ni\tau\longrightarrow
((\faz S)^* |S|^{i\overline{\tau}}y|\tilde{f}(|R|,\faz R)|S|^{1+i\tau}x)\in\C
\]
is continuous and holomorphic  on  $\Sigma$.

Consider for a while case $\tau\in\R$. Since $(R,S)\in D_{\cH}$, 
it follows that 
\[
((\faz S)^* |S|^{i\overline{\tau}}y|\tilde{f}(|R|,\faz R)|S|^{1+i\tau}x)=
((\faz S)^* y|S^{-i\tau}\tilde{f}(|R|,\faz R)|S|^{1+i\tau})=\]
\[=(y|\tilde{f}(e^{\hbar\tau}|R|,e^{\hbar\tau}\faz R)Sx).\]
The above equality remains true after  analytic continuation to the 
 $ \tau\in \ov{\Sigma}$.
In particular for  $\tau=i$ we get
\[
(S^*y|f(R)x)=
(y|\tilde{f}(e^{i\hbar}|R|,\faz R)Sx)=(y|f(q\cdot q R)Sx).
\]
The above formula holds for all  $y\in D(S^*)$. 
Moreover, 
$f(R)x\in D(S)$ and  $Sf(R)x=f(q\cdot qR)Sx$.\hfill\qed

For any  $z\in\Gamma$ 
 we define
\beq
\label{l}
\ell (z)=\log|z|\ .
\eeq
Let us define a new space of  holomorphic functions  
\beq
\label{Hzwykle}
H=\left\{
f\in C(\Gamma):
\mbox{for any  }\lambda>0\mbox{  function }e^{-\lambda\ell(r)^2}f(r)\in
 H_{bounded}
\right\}
\ .
\eeq
Note that  $H$ is a  vector space. Moreover, if we define 
  for any  $f\in H$ and $z\in\Gamma$  
\[f^*(z)=\ov{f(q\cdot q\ov z)}\ ,\]
then $f^*\in H$ and $(f^{*})^*=f$.
By the above consideration and functional calculus for normal operators,
 for any operator $R$ with spectrum contained in  $\ov\Gamma$, 
 we get 
\beq
\label{inwodopR}
f^*(R)=f(q\cdot q R^*)^*\ ,
\eeq
for any function  $f\in H$. 

Observe that for any  $\lambda>0$, 
 function 
$e^{-\lambda\ell(\cdot)^2}\in H_{bounded}$.
 
Let $(R,S)\in D_{\cH}$. Then for any   $\lambda>0$,  operator 
$e^{-\lambda\ell(R)^2}$ is 
 bounded  and converges strongly to  $I$ when  
 $\lambda\rightarrow+0$. 
So the set
\beq
\label{Do}
D_0=\bigcup_{\lambda>0} e^{-\lambda\ell(R)^2}D(S),
\eeq
 is dense in  $\cH$, because  $D(S)$ is dense in $\cH$.
 Moreover, we have
\bst
\label{QE3.2}
Let $f\in H$. Then
\bit
\item[0.] $D_0\subset D(f(R))$,

\item[1.] $f(R)D_0\subset D_0$, 

\item[2.] $D_0\subset D(S)$,

\item[3.] $SD_0\subset D(f(q\cdot qR))$.
\eit
\est

\dow 
Let  $\lambda>0$. One can easily check that function   
 $g(r)=e^{-\lambda\ell(r)^2}$ satisfies assumptions of 
 Proposition \ref{QE3.1}. Therefore  
 $e^{-\lambda\ell(R)^2}D(S)\subset D(S)$ and we thus proved 
 point 2.

Let now  $f\in H$ and  $\lambda>0$. Then the function 
\beq
\label{g2}
g(r)=f(r)e^{-\lambda\ell(r)^2}
\eeq 
is bounded. Therefore  
 $e^{-\lambda\ell(R)^2}D(S)\subset
D(f(R))$ and point 0. follows. 

One can easily check that function $g(r)$ defined by \rf{g2} 
 satisfies assumptions of 
 Proposition \ref{QE3.1}. Therefore 
\[
f(R)e^{-\lambda\ell(R)^2}D(S)\subset D(S).
\]
Moreover
\[
f(R)e^{-2\lambda\ell(R)^2}D(S)\subset e^{-\lambda\ell(R)^2}D(S)\subset D_0
\]
 and point 1. follows.

Let  $x\in D_0$. Then
\beq
\label{xprime}
x=e^{-\lambda\ell(R)^2}x',
\eeq
where  $x'\in D(S)$ and  $\lambda>0$. By  \rf{QE3.1wz} 
\[
e^{-\lambda\ell(q\cdot qR)^2}Sx'=Sx.
\]
If  $f\in H$ then function $f(r)e^{-\lambda\ell(r)^2}$ is bounded,
$f(q\cdot qR)e^{-\lambda\ell(q\cdot qR)^2}\in B(\cH)$ and 
$Sx=e^{-\lambda\ell(q\cdot qR)^2}Sx'\in D(f(q\cdot qR))$. 
This completes the proof of point 3.\hfill\qed

Let  $(R,S)\in D_{\cH}$ and let  $f\in H$. The just proved Proposition 
 shows that operators  $S\circ f(R)$ 
 and  $f(q\cdot q R)\circ  S$ are densely defined, in particular 
  their domains contain   $D_0$. 
We will prove that these  operators are closable. To this end we show that 
 their adjoints are densely defined. It means that adjoints of these adjoints 
 are well defined  - and they are  exactly closures of the considered 
operators.

Using formula   \rf{inwodopR} we deduce that   
\[(S\circ f(R))^*\supset
f^*(q\cdot qR^*)\circ S^*\mi  
(f(q\cdot qR)\circ S)^*\supset S^*\circ f^*(R^*)\ .\]
Moreover, Proposition \ref{p2.2} yields that if $(R,S)\in D_{\cH}$, 
then also $(R^*,S^*)\in D_{\cH}$. We know also that if  $f\in H$, 
then $f^*\in H$, too.
Therefore, by the Proposition we just proved,
operators $(S\circ f(R))^*$ and $(f(q\cdot qR)\circ S)^*$ 
 are densely defined . Hence   $S\circ f(R)$ and $f(q\cdot q R)\circ
 S$ are closable operators. 
We will denote their closures 
 by $Sf(R)$ and $f(q\cdot qR)S$, respectively.
\btw
\label{QE3.3}
Let  $(R,S)\in D_{\cH}$ and let  $D_0$ be defined by
 \rf{Do}. Then for any function $f\in H$
we have
\bit
\item[0. ] $D_0$ is a  core  
 for $f(q\cdot qR)S$,
\item[1. ] $\left(f(q\cdot qR)S\right)^*=S^*f^*(R^*)$,
\item[2. ] $f(q\cdot qR)S\subset Sf(R)$,
\eit
\etw
\dow  
{\bf Ad 1. } We already know that  $(f(q\cdot qR)S)^*\supset S^*f^*(R^*)$,  
 it is enough to prove that the opposite inclusion holds.
Let  $y\in D((f(q\cdot qR)S)^*)$ and
\mbox{$z=(f(q\cdot qR)S)^*y$.} Then 
\beq
\label{Ad1i0} 
(y|f(q\cdot qR)Sx)=(z|x)
\eeq
for any  $x\in D(f(q\cdot qR)S)$. In particular
 (see points  2. and 3. in Proposition
\ref{QE3.2}) the above relation holds for all $x\in D_0$. 
Therefore
\[
(y|f(q\cdot qR)Se^{-\lambda\ell(R)^2}x')=(z|e^{-\lambda\ell(R)^2}x')
\]
for any  $\lambda>0$ and $x'\in D(S)$. Using  \rf{QE3.1wz} we get 
\[
(y|f(q\cdot qR)e^{-\lambda\ell(q\cdot qR)^2}Sx')=
(z|e^{-\lambda\ell(R)^2}x')
\ .\]
By \rf{inwodopR}
\[
(f^*(R^*)e^{-\lambda\ell^*(R^*)^2}y|Sx')=
(e^{-\lambda(\ell(R)^*)^2}z|x').
\]
This relation holds for for any $x'\in D(S)$. 
Hence we obtain 
\[f^*(R^*)e^{-\lambda\ell^*(R^*)^2}y\in D(S^*)\mi 
e^{-\lambda\ell^*(R^*)^2}y\in D(S^*f^*(R^*))\ .\]
Moreover
\[
S^*f^*(R^*)e^{-\lambda\ell^*(R^*)^2}y=e^{-\lambda(\ell(R)^*)^2}z.
\]
The above formula holds for any  $\lambda>0$ and  
 the operator $S^*f^*(R^*)$ is closed. When  $\lambda\rightarrow
+0$ we obtain  $y\in D(S^*f^*(R^*))$  and  $S^*f^*(R^*)y=z$. 
Thus we proved that 
$(f(q\cdot q R)S)^*\subset S^*f^*(R^*)$ and hence point 1. follows. 

{\bf Ad 0. } Observe that by   \rf{Ad1i0} one can restrict
 $x$ to  $D_o$.  It shows that  
  $f(q\cdot qR)S$ and its restriction to  $D_o$ have the same adjoints. 
Therefore they have the same closures, too. In other words, $D_o$ is
 a core for  $f(q\cdot qR)S$.

{\bf Ad 2.} Let  $x\in D_o$. Then  $x=e^{-\lambda\ell(R)^2}x'$,
 where  $x'\in D(S)$ and $\lambda>0$. Formula  \rf{QE3.1wz} yields that  
$Sf(R)x=f(q\cdot qR)e^{-\lambda\ell(q\cdot qR)^2}Sx'$.
In particular for  $f=1$ we obtain $Sx=e^{-\lambda\ell
(q\cdot qR)^2}Sx'$.
Comparing these two formulas we get  
\[
f(q\cdot qR)Sx=Sf(R)x
\]
This formula holds for any  $x\in D_o$. Remembering that   $D_o$ 
is a core for  $f(q\cdot qR)S$, 
 we obtain $f(q\cdot qR)S\subset Sf(R)$.\hfill\qed

Before we proceed to discuss   properties of pairs  $(R,S)\in D_{\cH}$,
we give an example of such operators.
\bpr[The most important one: Schr\"odinger's pair]
\label{paraSch}
Let $\cH=
\Lk (\Gamma,d\gamma)$.
Then for any  $z\in\Gamma$, $z=q^kr$,
\[ (Rf)(z)=z f(z),\]
and
\[ (\faz R f)(q^kr)=q^k  f(q^kr)\mi (|R|f)(q^kr)=r f(q^k r) ,\]
where
\[D(R)=\{f\in\lk(\Gamma)\;:\;\sum_{k=1}^N\int_0^\infty r |f(q^kr)|^2dr
<\infty\} .\] 
Moreover
\[ (S f)(z)= f(q^{-1}\cdot q^{-1} z),\]
and
\[ (\faz S f)(q^{k} r)=  f(q^{k-1}r)\mi (|S|f)(q^k r)=f(q^{-1}\cdot  q^k r) .\]
The domain of $ D(|S|)$ consists of all functions  
 $f\in\Lk(\Gamma)$, such that there exists a  function $g\in \Lk (\Gamma)$ 
and a   function $\tilde{f}$ holomorphic in $\Lambda_o\times q^{\Z}$ 
 and such that for any 
  $k=0,1,...,N-1$ we have 
\[\lim_{\fil\rightarrow 0^{-}} \tilde{f}(e^{i\fil}r,q^k)=f(q^kr)\mi
 \lim_{\fil\rightarrow -\frac{2\pi}{N}^{+}} \tilde{f_k}(e^{i\fil}r,q^k)=
g(q^{k-1}r)\ ,\]
where limits are taken in  $\Lk$-norm.
Moreover, for fixed $\fil$ not equal multiply of
$\frac{2\pi }{N}$ and for any $\lambda>0$, function $e^{-\lambda 
\ell(re^{i\fil})^2}\tilde{f}(re^{i\fil})$ should be bounded.
From now on we will use notation
\[f(q^{-1}\cdot q^kr)=g(q^{k-1}r)\ .\]
One can prove that thus  defined operators $R$ and  $S$ satisfy  
Definition \ref{Dmu}, i.e. $(R,S)\in D_{\Lk(\Gamma)}$.
\epr

The example above is crucial and all other examples are 
 built up from this one.
More precisely, by  Stone-von Neumann Theorem, 
  every pair $(R,S)\in D_{\cH}$ \label{rpS} is a direct sum of a certain 
number of copies of pairs, which are all unitary equivalent to 
the  Schr\"odinger's pair.

We will use often the following remark and proposition.
\buw
\label{azu2}
Checking one by one conditions in  Definition 
\ref{Dmu} one can easily show that if  \mbox{$(R,S)\in D_{\cH}$}  
$\gamma_1,\;\gamma_2\in \Gamma$, then also 
 $(\gamma_1 R,\gamma_2 S)\in D_{\cH}$.
\euw

\bst[Proposition 2.2  \cite{az+b}]
\label{p2.2}
If $(R,S)\in D_{\cH}$ , then also
 \mbox{ $(R^*,S^*)\in D_{\cH}$} and $(S^{-1},R)\in D_{\cH}$ and 
$(S,R^{-1})\in D_{\cH}$ and $(R,SR)\in D_{\cH}$.
\est
\section{The quantum exponential function for the  "az+b" group}
\label{sec:theq}
As we explained in Introduction, a sum of closed  operators 
 may not be closed. Therefore, we will not consider $R+S$ itself but its 
appropriate closure.
Contrary to the case of selfadjoint $R$ and $S$ considered in 
\cite{qexp,paper1}, closure of the  sum $R+S$ is the desired operator,
 i.e. it is  normal and  its spectrum is contained in  
$\ov\Gamma$ (Theorem 2.4 w \cite{az+b}).
Let $R\dot{+}S$ denote closure of the sum $R+S$.
We give now a useful formula for $R\dot{+}S$, where 
$(R,S)\in D_{\cH}$,  involving 
 the quantum exponential function $F_N$. 
The definition of  $F_N$ is given below.
\bst[\cite{az+b}]
\label{r+s}
Let $(R,S)\in D_{\cH}$. Then
\[R\dot{+} S=F_N(S^{-1}R)^*SF_N(S^{-1}R)
=F_N(R^{-1}S)RF_N(R^{-1}S)\]
\est
The special function $F_N\;:\;\Gamma\;\rightarrow \; \C$
 is given by \cite{az+b}
\beq
\label{FN}
F_N (q^k r)=\left\{
\ba{cc}
\prod_{s=1}^{\frac{k}{2}}\left(\frac{1+q^{2s}r}{1+q^{-2s}r}\right)
\frac{f_o(qr)}{1+r}& {\rm \ for  \ } 2|k\\
\prod_{s=0}^{\frac{k-1}{2}}\left(\frac{1+q^{2s+1}r}{1+q^{-2s-1}r}\right)
f_o(r)& \;{\rm \ for  }\; 2\not{|}k
\ea
\right.\ ,
\eeq
where
\beq
\label{deffo}
f_o(z)=\exp \left\{\frac{1}{\pi i}\int_0^{\infty}\log(1+a^{-\frac{N}{2}})
\frac{da}{a+z^{-1}}\right\}\ ,
\eeq
for any $z\in \C\setminus\{\R_-\cup \{0\}\}$.

Note that if  $N$ would be an odd number,  this definition would not be  
 good - values of function $F_N$ at the point $q^kr=q^{k+N}r$ 
 could be calculated in two different ways giving different results. 
Whereas for   $N$ even we get, regardless of the way of calculating it,  
 all the time the same value at the same point, i.e. $F_N$ 
is well defined for an even $N$.

By \rf{deffo} it follows that relation between  $f_o$ and 
 the special function 
 $\vt$ used in \cite{qexp,paper1} is given by
\beq
\label{zwiazek}
f_o(q^kr)=
V_{\frac{N}{2}}(\log r+i\hb k)^2\ , 
\eeq
where $r\in\R_+$ and $k=0,1,\dots ,N-1$.

From  \cite{qexp} we know that  function $V_{\frac{N}{2}}$ is holomorphic  in 
the stripe $\Im x <\pi$, so in particular it is continuously differentiated 
  along  lines  $x=\log r$ and $x=\log r+i\hb$, for any  $r\in\R_+$. 
 
Moreover, by  (1.37) and  (1.38)  \cite{qexp}
 it follows that
\beq
\label{rozwVt}
V_{\frac{N}{2}}(\log r+il\hb)=1+\frac{q^lr}{2i\sin \hb}+{\cal R}_o(q^lr)
\ ,
\eeq
where
\[\lim_{r\rightarrow 0}\frac{{\cal R}_o(q^lr)}{r}=0\ \]
and  $l=0$ or $l=1$.

Hence
\[\lim_{r\rightarrow 0}V_{\frac{N}{2}}(\log r)^2=1=
\lim_{r\rightarrow 0}V_{\frac{N}{2}}(\log r+i\hb)^2\ .\]
Therefore setting $F_N (0)=1$ 
 will make  $F_N$ continuous on  
 $\ov\Gamma$.

Moreover, by (1.7)  \cite{az+b}
\[|F_N(\gamma)|=1 \dla \gamma\in\ov\Gamma\ .\]
In order to prove the theorem of this Section we will need  
 a formula for an expansion of  $F_N$ around $0$.
Let us derive a formula for  derivative at  the point zero for $k$ even
\[F^\prime_N(q^kr)=\left(\frac{1+q^2r}{1+q^{-2}r}\right)^\prime
\left(\frac{1+q^4r}{1+q^{-4}r}\right)\dots 
\left(\frac{1+q^kr}{1+q^{-k}r}\right)\frac{f_o(qr)}{1+r}
+\dots\]
\[\dots+\left(\frac{1+q^2r}{1+q^{-2}r}\right)
\left(\frac{1+q^4r}{1+q^{-4}r}\right)\dots 
\left(\frac{1+q^kr}{1+q^{-k}r}\right)^\prime\frac{f_o(qr)}{1+r}+
\prod_{s=1}^{\frac{k}{2}}\left(\frac{1+q^{2s}r}{1+q^{-2s}r}\right)
\left(\frac{f_o(qr)}{1+r}\right)^\prime
\]
Moreover
\[\left(\frac{1+q^{2s}r}{1+q^{-2s}r}\right)^\prime=
\frac{2i\Im q^{2s}}{(1+q^{-2s}r)^2} \mi 
\left(\frac{f_o(qr)}{1+r}\right)^\prime=
\frac{(1+r)f^\prime_o(qr)-f_o(qr)}{(1+r)^2}\ .\]
By \rf{zwiazek} and \rf{rozwVt} one can calculate  right 
 derivatives of  function  $f_o$ at the point zero
\[
f_o^\prime(r)|_{r=0}=\frac{1}{i\sin\hb }
\mi
f_o^\prime(qr)|_{r=0}=\frac{q}{i\sin\hb }\]
Hence
\[
F^\prime_N(q^kr)|_{r=0}=\frac{-iq}{\sin\hb }-1+2i\sum_{s=0}
^{\frac{k-1}{2}}\Im q^{2s+1}
=\frac{\cos \hb}{i\sin\hb }+2i\sum_{s=1}^{\frac{k}{2}}\Im q^{2s}
\]
Let us compute
\[\sum_{s=1}^{\frac{k}{2}}\Im q^{2s}=\Im q^{2}\sum_{s=0}^{\frac{k}{2}-1}q^{2s}
=\Im\left( q^{2}\frac{1-q^k}{1-q^2}\right)=
\Im\left( q\frac{1-q^k}{q^{-1}-q}\right)=\]
\[=\left( \frac{\Im(i(q-q^{k+1}))}{2\sin\hb}
\right)
=\frac{\cos\hb-\cos (k+1)\hb}{2\sin\hb}\]
Hence
\[F^\prime_N(q^kr)|_{r=0}=\frac{\cos\hb}{i\sin\hb}+
\frac{-\cos\hb+\cos (k+1)\hb}{i\sin\hb}=\frac{\cos (k+1)\hb}{i\sin\hb}\]
Finely
\[F^\prime_N(q^kr)|_{r=0}=
\frac{q^{k+1}+q^{-k-1}}{2i\sin\hb}\ .\]
 Similar calculations show that the above formula remains true also
 for  $k$ odd.
From the Taylor formula we obtain  
 an expansion of  $F_N$ around  $0$
for $\lambda \in \R_+$ and  $t\in \Gamma$
\beq
\label{ap}
F_N(\lambda t)=1+\frac{\lambda}{2i\sin \hb}(q t+\ov q 
 \ov t) +{ r}(\lambda t)\lambda |t|\ ,
\eeq
where
\beq
\lim_{\lambda\rightarrow 0} r(\lambda t)=0\ .
\label{gr0}
\eeq
Moreover, using  \rf{ap} one can deduce that
\[\lim_{\lambda\rightarrow +\infty} r(\lambda t)=
-F_N^\prime(\lambda t)|_{t=0}\ .\]
Obviously, function $r$ is  continuous on  $\gamma$. Hence, 
 there exist a constant $M$,independent of  $\lambda$, and such that
\beq
\label{ogrr}
|r(\lambda t)|<M, \mbox{\rm for any } t\in \Gamma.
\eeq
By \rf{ap} we have
\beq
\label{ap2}
\frac{F_N(\lambda t)-1}{\lambda}=
\frac{q  t+\ov q \ov t}{2i\sin\hb}+ { r}(\lambda t) |t|\ ,
\eeq
in particular
\beq
\label{grl0}
\lim_{\lambda\rightarrow 0}\frac{F_N(\lambda t)-1}{\lambda}=
\frac{q t+\ov q \ov t}{2i\sin\hb}\ .
\eeq
Our next objective is to derive  \rf{zal}, which we will use to 
 transform \rf{dorhs}.
Let $dE(t)$ be a  spectral measure   \cite[Chapter 8]{reedI} 
of normal operator $T$. 
According to rules of functional calculus  \cite[Theorem VIII.6]{reedI}, 
a continuous, bounded function of  an operator $T$ has form
\[\left(u|\;\left|\frac{F_N(\lambda T)-I}{\lambda}\;\right|v\right)=
\int_{\Gamma}\frac{F_N(\lambda t)-1}{\lambda}\left(u|\;dE(t)\;|v\right) \ , \]
for any $u\in\cH$ and $v\in D(T)$.

Moreover,  by \rf{ap2}
\beq
\label{dluugi}
\int_{\Gamma}\frac{F_N(\lambda t)-1}{\lambda}\left(u|\;dE(t)\;|v\right)=
\int_{\Gamma}
\frac{q  t+\ov q \ov t}{2i\sin\hb}\left(u|\;dE(t)\;|v\right)+ 
\int_{\Gamma}{ r}(\lambda t)  \left(u|\;|t|\;dE(t)\;|v\right)\ .
\eeq
Observe that  $\left(u|\;|t|dE(t)\;|v\right)=
\left(u|\;dE(t)|T|\;|v\right)$ and 
 the measure $\left(u|\;dE(t)|T|\;|v\right)$ is finite.

Constant function  $M$, which majorises  function $r$ (se \rf{ogrr}), 
is integrable with respect to the measure 
$ \left(u|\;|t|dE(t)\;|v\right)$ and  is the function  majorising  
$r (\lambda t)$.
 
Therefore, by the Lebesgue dominated convergence theorem 
\cite[Theorem I.16]{reedI} and by  \rf{gr0}
\[\lim_{\lambda\rightarrow 0}\int_{\Gamma}{ r}(\lambda t)  \left(u|\;|t|dE(t)\;|
v\right)=0\ .\]
Hence by \rf{dluugi}
\[\lim_{\lambda\rightarrow 0}\int_{\Gamma}\frac{F_N(\lambda t)-1}{\lambda}\left(u|\;dE(t)\;|v\right)=
\frac{1}{2i\sin\hb}\int_\Gamma(qt+ \ov{q}\ov t)\left(u\left|\;dE(t)\;\right|v\right)\ , \]
so
\[
\label{zal}
\lim_{\lambda\rightarrow 0}\left(u\left|\frac{F_N(\lambda T)-I}{\lambda}\right|v\right)
=\frac{1}{2i\sin\hb}\left(u\left|(q T+ \ov q  T^*)\right|v\right)\ ,
\]
for any $u\in\cH$ and $v\in D(T)$.
This formula will be used in proof of the main theorem of this Section.

The proposition below explains why $F_N$ is called 
 the quantum exponential function
\bst[Theorem 2.6, \cite{az+b}]
\label{ee}
Let $(R,S)\in D_{\cH}$. Then
\[
\label{eeaz}
F_N(R)F_N(S)=F_N(R\dot{+} S)
\]
\est
Later on in this Section we will prove that the quantum exponential 
 function is    
 unique (up to a  parameter) solution of  \rf{eeaz}. 
\section{Normal extensions of  $\mu R+RS$}
\label{sec:nor}
\bst
Let $(R,S)\in D_{\cH}$ and let $\mu\in \C$.
Operator $\mu R+RS$ has a  normal extension  if and only if, when   
 $\mu\in\Gamma$.
\label{nieu}
\est
\dow
$\Leftarrow$  Obvious, because according to Proposition  \rf{p2.2} we have 
$(RS,S)\in D_{\cH}$, and on account of  Remark \ref{azu2} 
it follows that $(RS,\mu S)\in D_{\cH}$, 
so
 $\mu S\dot{+}RS$ is by Theorem 2.4  \cite{az+b}   
 a  normal extension of   $\mu S+RS$
.\\
$\;\Rightarrow$ 
We show that for any  $\mu\in\C\setminus \Gamma$ operator $\mu S + RS$ 
 has not a  normal extension.
Let us introduce notation
\[Q=\mu S +  RS.\] 
We first prove that   operator 
 $Q$ is closed.

Let us define
\[m:=\sup_{z\in\Gamma}\frac{1}{|\mu+z|}= \frac{1}{d(-\mu,\Gamma)}\ ,\]
where $d(-\mu,\Gamma)$ denotes Euclidean distance  of the point  
$-\mu$ 
 from the set $\Gamma$.
Note that $m<\infty$, since 
 $\Gamma$ is closed and  $-{\mu}\notin\Gamma$.
Hence for any $z\in \Gamma$ we have 
\[c|\mu + z|\geq 1.\]
We know that  spectrum of operator $R$ is contained in $\ov\Gamma$ and 
that $R$ is  normal operator, i.e. there exists a representation 
in which it is a multiplication 
 by a $\ov\Gamma$-valued function operator.
Hence for any  $\psi\in D(\mu i +R)$
\[
c||(\mu i +R)\psi||\geq ||\psi||,
\]
 so for any $f\in D(Q)$ we have
\beq
\label{dq1}
c||Qf||\geq ||Sf||.
\eeq
 
Consider a sequence  $\phi_n$  of elements of $D(Q)$  norm-converging 
 to $\phi$ and  such that  $Q\phi_n$ is norm-converging to certain $y$. 
Since $D(Q)\subset D(S)$ and by inequality 
 \rf{dq1}, we conclude that  also $\phi\in D(S)$ and   $S\phi_n$ is 
norm-converging to $S\phi$, as operator $S$ is normal, so it is closed.
Since $\phi_n$ belongs to the domain of operator $Q$, hence it belongs  
 also to the domain of operator $RS$. Moreover,  $RS\phi_n$ 
is converging to $y-\mu S\phi$. 
As we know that operator $RS$ is closed we see that 
$\phi$ belongs  
 also to the domain of $Q$ 
and $y=RS\phi+ \mu S\phi=Q\phi$.
It means that  operator  $Q$ is closed.

By definition, the domain of  operator $Q^*$ consists of all  
   $x\in\cH$ such that  there exists  $w\in \cH$ such that for any  
 $y\in D(Q)$ we have
 \beq
\label{defdzQ*}
\langle x,Qy\rangle =\langle w,y\rangle\ .
\eeq
Set
\[f(r)=\mu+ q^{-2} r\ .\]
Observe that function $f$ belongs to space  $H$.
With this notation
\[Q=f(q\cdot q R)\circ S\ .\]
Since $D_o$ is a core  for $Q$,
 one can assume in formula  \rf{defdzQ*} that
$y\in D_o$. Every element  $y\in D_o$ has form 
$y=e^{-\lambda \ell (R)^2}z$ for certain $\lambda >0$ and $z\in D(S)$. 
Hence
\[\langle x,Qy\rangle=\langle x|f(q\cdot qR)S 
e^{-\lambda \ell (R)^2}z\rangle\]
Since function $e^{-\lambda \ell (\cdot)^2}\in H_{bounded}$, \rf{QE3.1wz}
 shows that 
 \[\langle x,Qy\rangle=\langle x|f(q\cdot qR) 
e^{-\lambda \ell (q\cdot qR)^2}S z\rangle\ .\]
Moreover, since operator $f(q\cdot qR) 
e^{-\lambda \ell (q\cdot qR)^2}$ is bounded, we see that the
 domain of its adjoint  operator  is the whole space $H$.
 Hence by \rf{QE3.1wz}
\[\langle x,Qy\rangle=\langle
f^*(R^*) e^{-\lambda \ell^* (R^*)^2} x| 
S z\rangle
=\langle
f^*(R^*)  x| 
e^{-\lambda \ell (q\cdot qR)^2}S z\rangle \ . \]
Consequently, if  vector $x\in \cH$ belongs to the domain $Q^*$, then
\beq
\label{warQ*}
(\ov\mu i +R^*)  e^{-\lambda \ell^* (R^*)^2} x\in D(S^*)=D(S)\ .
\eeq
Moreover
\beq
\label{warQ*2}
S^*(\ov\mu i +R^*)  e^{-\lambda \ell^* (R^*)^2} x=
e^{-\lambda (\ell (R)^*)^2}Q^* x\ .
\eeq
Let $$Q^{\prime}=(Q^* |_{D(Q)})^*\ . $$
On account of Proposition \rf{QE3.3} point 2 
\[Q^*=(f(q\cdot q R)\circ S)^*=S^* f^*(R^*)\ .\]
We will prove that
\beq
\label{wznaQ*ob}
Q^* |_{D(Q)}=f(q\cdot q R)^*\circ S^*\ . 
\eeq
Vector $x\in D(Q)$ belongs to the domain of  operator $Q^* |_{D(Q)}$, if  
 for any  $y\in D(Q)$ there exists  $w \in H$ such that
\[ \langle x|Qy\rangle =\langle w|y\rangle \ . \]
We know, that $D_o$ is a core for $Q$.
Therefore we can set $y\in D_o$ in the above formula. 
 Every element $y\in D_o$ has form 
\[y=e^{-\lambda \ell (R)^2}z\ ,\]
 for any $\lambda >0$ and $z\in D(S)$. Hence
\[\langle x,Qy\rangle=\langle x|f(q\cdot q R)S 
e^{-\lambda \ell (R)^2}z\rangle\]
Since $z\in D(S)$ and   $e^{-\lambda \ell (\cdot )^2}
\in H_{bounded}$, formula \rf{QE3.1wz}
 implies
\[ \langle x| f(q\cdot q R)e^{-\lambda \ell (q\cdot qR)^2} S z\rangle
\ .\]
Moreover function $f(\cdot)e^{-\lambda \ell (\cdot )^2}
\in H_{bounded}$, so using once more   \rf{QE3.1wz} 
we obtain
\[ \langle x| Sf( R)e^{-\lambda \ell (R)^2}  z\rangle
\ .\]
Since $x\in D(Q)$, it follows that  $x\in D(S)=D(S^*)$. Hence 
\[ \langle x| Sf( R)e^{-\lambda \ell (R)^2}  z\rangle
= \langle S^* x| f( R)e^{-\lambda \ell (R)^2}  z\rangle\ .\]
Moreover $Sx\in D(f(q\cdot qR))=D(f(R)^*)$. It is easily seen that 
 the domain of  $f(R)^*$ is invariant  with respect to the action 
 of  operator
$\faz S$. Hence also $S^*x\in D(f(R)^*)$.
Consequently
\[\langle S^* x| f( R)e^{-\lambda \ell (R)^2}  z\rangle
=\langle  f( R)^*S^* x|e^{-\lambda \ell (R)^2}z\rangle=
\langle  f( R)^*S^* x|y\rangle\ .\]
Thus we proved \rf{wznaQ*ob}.

We proceed to derive condition which  $x\in \cH$ has to satisfy to belong 
 to the domain of  $Q^\prime$.
By definition,  $x\in D(Q^\prime)$ if for any $y\in D(Q)$
\[ \langle x|Q^*y\rangle =\langle Q^\prime x|y\rangle \ . \]
By \rf{wznaQ*ob} and \rf{inwodopR}
\[ \langle x|Q^*y\rangle =\langle x|f^*(R^*)\circ S^*y\rangle  \ . \]
On account of  \rf{QE3.3} point 1,  $D_o$ is a core of 
operator $Q^*|_{D(Q)}$. 
Hence we may assume that  $y\in  D_o$. 
Every element $y\in D_o$ has form  
\[y=e^{-\lambda \ell (R^*)^2}z\ ,\]
 for any $\lambda >0$ and $z\in D(S^*)$. Hence
\[\langle x,Q^*y\rangle=\langle x|f^*( R^*)S^* 
e^{-\lambda \ell (R^*)^2}z\rangle\ .\]
By \rf{QE3.1wz}, as previously, we obtain
\[\langle x,Q^*y\rangle=\langle x|f^*( R^*) 
e^{-\lambda \ell (q\cdot qR^*)^2}S^*z\rangle\ .\]
Moreover, since operator $f^*(R^*) 
e^{-\lambda \ell (q\cdot qR^*)^2}$ is bounded, the domain 
 of its adjoint is the whole Hilbert space $H$.

Hence by \rf{QE3.1wz}
\[\langle x,Q^*y\rangle=\langle
 e^{-\lambda \ell^* (R)^2}f(q\cdot q R) 
x| S^*z\rangle\ .\]
From this we conclude that if  $x\in D(Q^\prime)$, then   
\beq
(\mu I+R)e^{-\lambda \ell^* (R)^2}x\in D(S)\ .
\label{warQ'}
\eeq
Moreover
\beq
S(\mu I+R)e^{-\lambda \ell^* (R)^2}x=
e^{-\lambda (\ell (R^*)^*)^2}Q^{\prime}x
\label{warQ'2}\ .
\eeq
Observe also that by \rf{l} for any  $z\in\Gamma$
\[\ell^*(\ov z)=\ell^*(z)\ .\]
It implies that 
\[e^{-\lambda\ell^*(R^*)^2}=e^{-\lambda \ell^*( R)^2}\ ,\]
for any $\lambda >0$.
This result allows one to replace \rf{warQ*} with a more convenient,
 as we will soon see, condition
\beq
(\ov\mu I+R^*)e^{-\lambda \ell^* (R)^2}x\in D(S)\ .
\label{warQ*l}
\eeq
Assume that  $x\in D(Q^*)\cap D(Q^{\prime})$. 
It means that  $x$ satisfies simultaneously  \rf{warQ*l} and \rf{warQ'},
i.e.
\[(\ov\mu i +R^*)  e^{-\lambda \ell^* (R)^2} x\in D(S)
\mi  (\mu I+R)e^{-\lambda \ell^* (R)^2}x\in D(S)\ .\]
Because the domain of  operator $S$ is invariant with  respect to 
the action of the phase of $R$, the second condition yields that 
\beq
(\ov\mu (\faz R)^2  +R^*)  e^{-\lambda \ell^* (R)^2} x\in D(S)
\label{fazr}
\ .
\eeq
Since a domain of any linear operator  is a linear subspace of $\cH$, 
 hence  by \rf{warQ*l} and \rf{fazr} it follows that
\beq
(\ov\mu (\faz R)^2 -\mu I)  e^{-\lambda \ell^* (R)^2} x\in D(S)
\label{fazr2}
\ .
\eeq
Let $\rho$ be the  $N$th primitive root of unity, i.e. let 
$\rho \in \C$ and  $\rho^N=1$.
Let function $h$ be defined on the set of  $N$Th roots of unity 
and be given by
\[h(\rho)=\ov\mu \rho^2-\mu\ .\]
Note that $h(\rho)\neq 0$, since in order  
 to have equality $\ov\mu q^{2k}=\mu$, the  phase of  $\mu$ 
 would have to be equal to the multiply of $\frac{2\pi}{N}$.
But this contradicts our assumption that $\mu\not\in \Gamma$.

We know that $\faz R$ is a unitary operator which eigenvalues 
 are all  $N$th roots of unity.
Consequently the operator
\[h(\faz R)=\ov\mu (\faz R)^2-\mu\ \]
is invertible.
Moreover it is easily seen that the operator $h(\faz R)^{-1}$
 is given by
\[h(\faz R)=a_0 i + a_1 (\faz R) + a_2 (\faz R)^2 +\dots +
a_{N-1}(\faz R)^{N-1}  ,\]
where  $a_0, a_1, a_2, \dots, a_{N-1}\in \C$.

Therefore, because 
\[ h(\faz R) e^{-\lambda \ell^* (R)^2} x\in D(S)
\]
and the domain of the  operator $S$ is invariant with 
respect to the action of the 
phase of  operator $R$, it follows that
 \[h(\faz R)^{-1} h(\faz R) e^{-\lambda \ell^* (R)^2} x\in D(S)\ ,
\]
 and finally
\beq
\label{fazr3}
e^{-\lambda \ell^* (R)^2} x\in D(S)\ 
\eeq
for any $\lambda >0$.

Let  $\lambda=\lambda_1+\lambda_2$, where $\lambda_1,
\lambda_2\in\R_+$. 
By\rf{QE3.1wz}
\[ Se^{-\lambda \ell^* (R)^2}x=e^{-\lambda_1 \ell^* (q\cdot qR)^2}
Se^{-\lambda_2 \ell^* (R)^2}x\ .\]
Moreover the operator $Re^{-\lambda_1 \ell^* (q\cdot qR)^2}$ is bounded, 
hence for any $y\in\cH$ 
 $$ e^{-\lambda_1 \ell^* (q\cdot qR)^2}y\in D(R).$$
Moreover 
$$(\mu+\cdot )e^{-\lambda_1 \ell^* (\cdot )^2}\in 
H_{bounded}.$$ 
Hence by \rf{QE3.1wz} for any  $x\in\cH$ we have
\beq 
S(\mu I+R)e^{-\lambda \ell^* (R)^2}x=(\mu I+q^2R)
e^{-\lambda_1 \ell^* (q\cdot qR)^2}Se^{-\lambda_2 \ell^* (R)^2}x=
(\mu I+q^2R)Se^{-\lambda \ell^* (R)^2}x
\label{laczac1}\ .
\eeq
On the other hand by \rf{warQ'2}
\beq
S(\mu I+R)e^{-\lambda \ell^* (R)^2}x=e^{-\lambda (\ell (R^*)^*
)^2}
Q^{\prime}x
\label{laczac2}
\ .
\eeq
Combining \rf{laczac1} and \rf{laczac2} we obtain
\beq
(\mu I+q^2R)Se^{-\lambda \ell^* (R)^2}x
=e^{-\lambda (\ell (R^*)^*
)^2}
Q^{\prime}x
\label{laczac3}
\ .
\eeq
It has been already proved that the 
operator $Q=(\mu I+R)\circ S$ is closed.
Therefore on account of remark \ref{azu2} the operator $(\mu I+q^2R)\circ S$
 is also closed.
Hence    setting  $\lambda\rightarrow 0^+$ 
 into formula  \rf{laczac3}
 yields
\[x\in D((\mu I+q^2R)\circ S)=
D((\mu I+R)\circ S)=D(Q)\ .\]
We have have proved that  if  $x\in D(Q^*)\cap D(Q^{\prime})$, 
then $x\in D(Q)$.
Rephrasing, we have proved that
\[
D(Q)\supset D(Q^*)\cap D(Q^{\prime})\ . 
\]
On the other hand, from the definitions of  operators 
$Q^*$ and $Q^{\prime}$ it is obvious that  
\[
D(Q)\subset D(Q^*)\cap D(Q^{\prime})\ . 
\]
Consequently, we have proved that
\beq
\label{dq}
D(Q)= D(Q^*)\cap D(Q^{\prime})\ . 
\eeq
Note that the domain of   $Q^*$ is strictly greater than this of $Q$. 

In order to give an example of a function  
$\psi\in\Lk (\Gamma)$, such that 
\[\psi\in D(Q^*)\setminus D(Q)\ ,\]
let us assume that  $(R^*,S^*)$ is the Schr\"odinger pair 
(see Example  \ref{paraSch}).
Then operator $Q^*$ acts on functions $\psi$ from its domain as follows
\[ (Q^*\psi)(z)=(S^*((\ov\mu +z )\psi))(z)
\ .\]
Let function $\psi$ be for any $z\in \Gamma$ given by
\[\psi(z)=\frac{e^{-\ell (z)^2}}
{ z+\ov\mu}\ .\]
Observe that $\psi\in\Lk(\Gamma)$, since $\mu\not\in\Gamma$.
Moreover the function $\widetilde{(\ov\mu+\cdot)\psi}=e^{-\ell (\cdot)^2}$
 corresponding to the function   $(\ov\mu+\cdot)\psi$ in the   sense 
 described in Example \ref{paraSch},
is holomorphic and bounded inside $\Lambda_o\times q^{\Z}$.
Moreover, one can easily check  that this function satisfies also 
 other conditions  of belonging to the domain of  $|S|$ 
 given in Example  \ref{paraSch}.
Moreover,  $Q ^*\psi\in\Lk(\Gamma)$.
Consequently,  $\psi$ belongs to the domain of $Q^*$. 

However, $\psi$ does not belong 
 to the domain of  operator $Q$, since it does not  belong 
to the domain of  $|S|$. Function $\psi$ 
  does not  belong 
to the domain of  $|S|$, since its corresponding 
 ( in the sense explained  in Example \ref{paraSch}) 
function $\tilde{\psi}$ is 
  meromorphic (but it is not  holomorphic). 
Namely, function $\tilde{\psi}$ has a simple pole at the point  $z= -\ov\mu$.
It shows that $D(Q)_{\;\not}\!\!\subseteq D(Q^*)$ and $Q$ is 
closed, so $Q$ is not normal.

We repeat now  S.L. Woronowicz's  reasoning from the proof
  of  Theorem 2.2 \cite{oe2}.
assume, that $\tilde{Q}$ is a normal extension of $Q$.
Then 
\beq
\label{dq2}
D(Q)\subset D(\tilde{Q})=D(\tilde{Q}^*)\subset D(Q^*)\ .
\eeq
On the other hand
\[Q^* |_{D(Q)}\subset \tilde{Q}^*\subset Q^*\ ,\]
hence
\beq
\label{dq3}
D(\tilde{Q})\subset D(Q^{\prime}),
\eeq
since   $Q^{\prime}=(Q^* |_{D(Q)})^*$.

Using \rf{dq},\rf{dq2} and \rf{dq3} we obtain 
 that 
$$D(\tilde{Q})\subset D(Q).$$
But by assumption  $\tilde{Q}$ is a normal extension of $Q$, 
so 
$\tilde{Q}=Q$.

This last result implies that if  $Q$ is not normal,
it has not a normal extension. We recall that 
 we have already proved that $Q$  is normal.

We have thus proved that if  $\mu\notin\Gamma$ then the  operator 
$\mu S +RS$ has no  normal extension, which completes the proof.
\hfill\qed

\section{Solutions of the scalar exponential equation}
\label{sec:sols}
We proceed to the proof of the main theorem of this Section.
\btw
\label{hip}
Let $(R,S)\in D_{\cH}$ and let  $\;f : \ov\Gamma\;\rightarrow \;S^1\;$
 be a Borel function.  
The following conditions are equivalent
\bit
\item[1).]
\beq
\label{ees}
f(R)f(S)=f(R\dot{+} S)\ 
\eeq
\item[2).]
\beq
\label{fs}
f(z)=F_N(\gamma z)
\eeq
for any $\gamma\in\ov\Gamma$ and almost all $z\in \ov\Gamma$.
\eit
\etw
\dow 1). $\Leftarrow$ 2). By Remark \ref{azu2} we know  
$(\gamma R,\gamma S)\in D_{\cH}$ for any  $\gamma\in\ov\Gamma$. 
Moreover $ \gamma R\dot{+}\gamma S=\gamma (R\dot{+}S)$. 
Hence by Proposition \ref{ee} 
\[F_N(\gamma R)F_N(\gamma S)=F_N(\gamma (R\dot{+} S))\ ,\]
which shows that function \rf{fs} satisfies 1)..\\

2). $\Rightarrow$ 1). Applying Proposition \ref{r+s}
 to the right-hand-side of \rf{ees} we obtain
\beq
\label{pow}
f(R)f(S)=F_N(S^{-1}R)^*f(S)F_N(S^{-1}R)\ .
\eeq
Let  $\lambda > 0$ be an arbitrary real positive number.
 Then by Remark \ref{azu2} it follows that $(\lambda R,S)\in D_{\cH}$. 
Substituting $\lambda R$ instead of  $R$ into   
 equation \rf{pow} and setting 
  $T=S^{-1}R$ we obtain
\[f(\lambda R)f(S)=F_N(\lambda T)^*
f(S)F_N(\lambda T)\ .\]
Equivalently
\[f(\lambda R)f(S)-f(S)=F_N(\lambda T)^*
f(S)F_N(\lambda T)-F_N(\lambda T)^*
f(S)+F_N(\lambda T)^*
f(S)-f(S)\ .\]
(Note the all the operators above are bounded, so adding and subtraction them do not change theirs domains. However, if the operators involved would be unbounded, one could get a false inequality, because such a procedure
could change domains of operators.)

Dividing both sides by  $\lambda>0$ we get
\beq
\frac{f(\lambda R)-I}{\lambda}f(S)=F_N(\lambda T)^*
f(S)\frac{F_N(\lambda T)-I}{\lambda}
+\frac{F_N(\lambda T)^*
-I}{\lambda}f(S)\ .
\label{dorhs}
\eeq
Let $RHS$ and $LHS$ denote  right and  and left hand 
side of 
 \rf{dorhs} and let $x,z\in D(T)$. Then
\beq
\label{porl}
\lim_{\lambda\rightarrow 0}\left(z\left| RHS\right| x\right)=
\frac{-i}{2\sin\hb}\left(z\left|f(S)(q T+ 
\ov q T^*)\; \right. x\right)
+\frac{i}{2\sin\hb}\left(\;(q T+ \ov q T^*)z\;\left|
f(S)\right. x\right) .
\eeq
By the reasoning above
\bwn
\label{pople}
For any  $x,z\in D(T)$ there exists the limit 
\[\lim_{\lambda\rightarrow ^+0}\left(z\left|\frac{f(\lambda R)-I}{\lambda}
\right|f(S)x\right)
\ .\]
\ewn
We will use the Lemma below to write this limit in the  
 form convenient for future computation
\ble[Lemma 7.2 z \cite{qexp}]
\label{lmu}
 Let $f$ be a bounded Borel function on 
  $\R_+$ and let $Y$ be a selfadjoint positive operator
 acting on   Hilbert space $\cH$. Let $\Xi$ be 
 the set of  pairs $(z,x)\in D(Y)\times \cH$ such that 
there exists the limit
\[\lim_{\lambda\rightarrow ^+0}\left(z\left|\frac{f(\lambda Y)-I}{\lambda}
\right|x\right)
\ .\] 
If there is a  pair $(y,u)\in \Xi$ such that  $(Yy|u)\neq 0$, 
then there is also  a constant $\mu_f$ such that
\[\lim_{\lambda\rightarrow ^+0}\left(z\left|\frac{f(\lambda Y)-I}{\lambda}
\right|x\right)
=\mu_f(Yz|x)\]
 for any $(z,x)\in \Xi$, and  $\mu_f$ depends only on $f$.
\ele

Note the the Lemma above can be applied to the  function $f$ mentioned in  
 Theorem \ref{hip} and satisfying condition \rf{ees} (because  this 
 function is by assumption Borel and bounded) and to the   operator 
$Y=|R|$, which is clearly selfadjoint and positive. 

Let $\cH_k$ denote an eigenspace of $\faz R$ 
 corresponding to the eigenvalue $q^k$. 
We have the following decomposition of the Hilbert space  $\cH$:
\[ \cH=\cH_1\oplus\cH_2\oplus\dots\oplus\cH_N \ .\]
The restriction of Operator $R$  to the space $\cH_k$ is $q^k |R|$.
(Matrix representation  of phases of operators $(R,S)\in D_{\cH}$ 
 is discussed in Section \ref{sect:az2} and given  by \rf{frmac} 
and \rf{fsmac})
By Lemma \ref{lmu} applied  independently 
in every space $\cH_k$  to the function $f(q^k .)$ and the operator  
$|R|$  and by Observation \ref{pople} 
 it follows that 
\beq
\label{poprz}
\lim_{\lambda\rightarrow 0}\left(z\left|
\frac{f(\lambda q^k|R|)-I}{\lambda}\right|y\right)
=\frac{-i}{2\sin\hb}\left(\ov{\mu(q^k)}|R| z\left|\right. 
y\right)\ ,
\eeq
where $z\in D(R)\cap D(T)$ and $y\in f(S)D(T)$ and
$\mu$ is a complex-valued function defined on the set of all 
  $N$th roots of unity.
Observe that the function $\mu$ is determined  uniquely  
 by complex numbers $b_1,b_2, ..., b_N$ such that 
\beq
\label{postacmu}
\mu(\rho)=b_1+b_2 \rho+b_3 \rho^2+...+b_{N} \rho^{N-1}\ ,
\eeq
where $\rho\in \C$ and $\rho^N=1$.
Note that$\rho\in\Gamma$.

To make further computations easier we adopt different $\mu$ in order 
to have  
$\frac{-i}{2\sin\hb}$ before parentheses in \rf{poprz}.

Substituting  $y=f(S)x$, where $x\in D(T)$, to \rf{poprz}
\beq
\label{porp}
\lim_{\lambda\rightarrow 0}\left(z\left|\frac{f(\lambda R)-I}{\lambda}
 \right| f(S) x\right)=\lim_{\lambda\rightarrow 0}\left(z\left|LHS 
\right|x\right)=
\frac{-i}{2\sin\hb}\left(\mu(\faz R)^*|R|z\left|\right.
f(S)x\right)\ .
\eeq
Comparing \rf{porp} and \rf{porl} we conclude that
\beq
\label{podwie}
\left(\mu(\faz R)^*|R|z\left|f(S)\right. x\right)=
\left(z\left|f(S)(qT+ \ov q T^*)\right. x\right)
-\left((q T+ \ov q T^*) z\left|f(S)\right. x\right)\  ,
\eeq
for any $z\in D(R)\cap D(T)$ and $x\in D(T)$.

From Remark \ref{azu2} it follows that  if operators 
$(R,S)\in D_{\cH}$ satisfy the above equation they satisfy also 
 $(\rho R,S)\in D_{\cH}$.
Let us substitute $\rho R$ instead of  $R$ to the equation  \rf{podwie} and  
write function  $\mu$ in the form \rf{postacmu}
\[\ b_1\left(z|f(S)x\right)+  \rho b_2 \left(R^*z|f(S)x\right)+
\rho^2 b_3 \left((\faz R)^{-2}|R|z|f(S)x\right)+\dots \]
\[\dots +\rho^{N-2} b_{N-1} 
\left((\faz R)^{2}|R|z|f(S)x\right)+\rho^{N-1} b_{N}
\left(Rz|f(S)x\right)=\]
\[=\rho q\left\{\left(z\left|f(S)Tx\right)-\left(T^*z|f(S)\right. x\right)
\right\}+
\rho^{N-1}q^{-1}\left\{\left(z|f(S) T^*x\right)-\left(Tz|f(S)x\right)\right\}
\ .`\]
Since  $\rho^k$ are linearly independent  for different  
$k\in\{0,1,2,\dots,N-1\}$, it follows that  by comparing coefficients 
 of the same powers of $\rho$ on both sides  of the above equation,
we obtain that  $b_k=0$ for $k\neq 2$ and $k\neq N$ 
and
\beq
\label{r1}
 b_{N}q\left( Rz\left|\right.f(S)x\right)
=\left(z|f(S)T^*x\right)-\left(Tz|f(S)x\right)\ 
\eeq
and
\[
 b_{2}q^{-1}\left( R^* z|f(S)x \right)
=\left(z|f(S) Tx\right)- \left(T^*z|f(S)x\right)\ .
\]
Setting  $\gamma=\ov{ b_{N}q}$ and rewriting  
\rf{r1} in a slightly different form we get
\[\left((\gamma R+T)z|f(S)x\right)=
\left(z|f(S)T^*x\right)\ ,\]
for any $z\in D(R)\cap D(T)$ and $x\in D(T)$.
Let us set $y=f(S)x$, i.e. $y\in f(S)D(T)$.
Then $D(T)=D(T^*)$, since $T$ is normal and 
$D(f(S))=\cH$, since $f(S)$ is bounded. So
\[\left((\gamma R+T)z|y\right)=
\left(z|f(S)T^*f(S)^*y\right)\ ,\]
for any $z\in D(R)\cap D(T)$ and $y\in D(f(S)T^*f(S)^*)$.
We see that the operator $ \gamma R+ T$ is 
contained in the operator adjoint to  $f(S)T^* f(S)^*$
\beq
\gamma R+ T\subset \left( f(S)T^* f(S)^*\right)^*=f(S)Tf(S)^*\ .
\label{zaw}
\eeq
The operator on the right-hand-side is unitarily equivalent to 
 a normal  operator, 
so it is  normal. It means that the  operator $\gamma R+ T$ has a 
 normal extension. Proposition \ref{nieu} implies that in such a 
 case $\gamma \in \Gamma$. Consequently, by Remark  \ref{azu2} 
 we have $(T,\gamma R)\in D_{\cH}$.
By Theorem 2.4  \cite{az+b} the closure of  $\gamma R+ T$, denoted  
by $\gamma R\dot{+} T$,  is normal for any $\gamma \in \Gamma$.
Moreover, Proposition 
\ref{ee} yields
\beq
\label{r3}
\gamma R \dot{+} T= F_N(\gamma S)T F_N(\gamma S)^*.
\eeq
Because by  Proposition  \ref{nieu}  $\gamma R+ T$ has a   normal extension 
 and by Theorem 2.4  \cite{az+b} this  normal extension of  $\gamma R+ T$ 
is  $\gamma R\dot{+} T$,  
 by \rf{zaw} and since normal operators do not have 
normal extensions, one may conclude that
\beq
\label{r5}
\gamma R \dot{+} T= f(\gamma S)T f(\gamma S)^*.
\eeq
Comparing  \rf{r5} and \rf{r3} we obtain
\[F_N(\gamma S)T F_N(\gamma S)^*=f(S) T f(S)^*\ .\]
Hence
\[T F_N(\gamma S)^*f(S)=F_N(\gamma S)^*f(S) T  ,\]
so $T$ commutes with a bounded operator $F_N(\gamma S)^*f(S)$. 
By  
spectral theorem for normal  operators \cite[Theorem 13.33]{rudaf}  
 it follows that   $F_N(\gamma S)^*f(S)$ commutes also with 
functions $|T|^{it}$, where $t\in \R$.
Hence
\[|T|^{it} F_N(\gamma S)^*f(S)|T|^{-it}=F_N( \gamma S)^*f( S)  \ ,\]
so
\beq
\label{nzol}
 F_N(\lambda\gamma S)^*f(\lambda S)=F_N(\gamma S)^*f( S) \ ,
\eeq
where $\lambda =e^{-\hb t}$.
Considerations preceding  Lemma \ref{lmu} imply that for any  $(T,R)
\in D_{\cH}$
\[\lim_{\lambda\rightarrow 0}\left(z|f(\lambda R)\xi\right)=
\left(z|\xi\right)\ ,\]
for any  $z\in D(T)$ and $\xi\in f(S)D(T)$. 
Since  $(T,S)\in D_{\cH}$, hence we obtain
\[\lim_{\lambda\rightarrow 0}
\left(\zeta|F_N(\lambda\gamma S)^*f(\lambda S)\xi\right)=
\left(\zeta|\xi\right)\ ,\]
where $\zeta\in F_N(\gamma S)^* D(T)$ and  
$\xi\in f(S^{-1}) D(T)$. Note that  $D(T)$ is linearly dense 
in $\cH$, as the  operator $T$ is densely defined. 
Similarly, we see that  
 $f(S^{-1})D(T)$ and $ F_N(\gamma S)^* D(T)$  are densely defined, 
 since $f(S^{-1})$ and $ F_N(\gamma S)$ are unitary.
Hence $F_N(\lambda\gamma S)^*f(\lambda S)$ converges weakly to $I$, 
 when $\lambda$ 
 goes to $0$. Since the right-hand-side of \rf{nzol} is independent on
 $\lambda$, it follows that
\[
 F_N(\lambda\gamma S)^*f(\lambda S)=I\ ,\]
so
\[f(S)=F_N(\gamma S)\ ,\]
where $\gamma \in \Gamma$, which completes the proof.\hfill\qed

\section{Solutions of the general exponential equation 
 for the quantum  "az+b" group}
\label{sect:az2}
We prove now the generalization of the Theorem \ref{hip} to 
the case of  an operator-valued function  $f$. 
Our earlier results from \cite{paper1,phd} will 
 make this proof much simpler. 

Let $f$ be a  function  defined on   $\Gamma$ and such, that  
 for any $z\in \Gamma$,  $f(z)$ 
 is a unitary operator acting on a Hilbert space $\cK$, i.e. 
$f(z)\in{\rm Unit}(\cK)$.
We will call  $f$ a Borel function iff for any 
$\fil,\psi\in\cK$ the function 
\[z\rightarrow (\fil|f(z)\psi)\]
 is Borel.
Then for any normal  operator $R$ with spectrum contained in  $\Gamma$
 we  define  function  $f(R)$ by
\[f(R)=\int_\Gamma f(z)\te dE_R(z),\]
where $dE_R$ is the spectral measure of the operator $R$ and $f(z)$ 
is a unitary operator acting on a Hilbert space  $\cK$.
\btw
\label{hip2}
Let $f$ be a  Borel function  defined on  $\Gamma$ with values in 
unitary operators acting on  $\cK$ 
 and let  $(R,S)\in D_{\cH}$. Then  
\[\left(
\begin{array}{c} 
f(R)f(S)=f(R\dot{+} S)
\end{array}
\right)
\Longleftrightarrow
\left(
\begin{array}{c} \mbox{ there exists  }\mbox{ an invertible}\\
\mbox{ normal operator  $M$, }\\
\mbox{ such that  }\spe M\subset \Gamma\\
\mbox{ and   }f(z)=F_N(M  z)\\
\mbox{ for a. a. } z\in \ov\Gamma.
\end{array}
\right).\]            
\etw
\dow $\Leftarrow  \;$ If $f(z)=F_N(Mz)$, then $f(R)=F_N(M\te R)$, where  
\[F_N(M\te R)=\int_\Gamma F_N(Mz)\te dE_R(z),\]
where $dE_R$ is the   spectral measure of $R$.
By assumption 
 $(R,S)\in D_{\cH}$ and $M$ is an invertible  normal  operator
 such that 
$\spe M\subset \Gamma$. Let us set 
 $R^\prime=M\te R$ and $S^\prime=M\te S$. Then  
$(R^\prime,S^\prime)\in D_{\cH\te\cK}$.
Therefore by  Theorem  \ref{hip} we get
\[F_N(M\te R)F_N(M\te S)=F_N(M\te (R\dot{+}S))\ .\]
 $\Rightarrow  \;$
We will use the same method as in the proof of Theorem 2.6 \cite{paper1},
 i.e.  we show that 
\[f(z)f(x)=f(x)f(z)\ ,\]
for any $z,x\in\Gamma$.
For the Reader's convenience we recall why it is enough to prove this.
Observe that 
 if  $\dim \cK=k<\infty$, then  
 from commutation of  unitary operators $f(z)$ and $f(x)$ follows that 
 there exists orthonormal basis in which these  operators 
 are  represented by  diagonal matrices for any  $x,y\in\Gamma$. 
Thus the problem reduces to  
 finding solutions of  $k$ scalar equations  
\[f_o(R)f_o(S)=f_o(R\dot{+}S)\ ,\]
where  $f_o$ is a complex-valued  function defined on $\Gamma$.
The above reasoning can be generalized to the  case of arbitrary 
many dimensional separable Hilbert space  $\cK$.  
This is so because  operators $f(q^kr)$ and  $f(q^ls)$ belong  to 
 a commutative *-subalgebra of $B(\cK)$. Therefore, by  spectral 
theorem and its consequences  \cite[Chapter X]{dunfII},  
operators $f(q^kr)$ and $f(q^ls)$ have  the same  spectral measure 
\[f(q^kr)=\int_0^{2\pi}f_o(q^kr,t)dE_{\cK}(t)\mi 
f(q^ls)=\int_0^{2\pi}f_o(q^ls,t)dE_{\cK}(t).\]
Hence
\[f(R)=\int_\Gamma \int_0^{2\pi}f_o(q^kr,t)dE_{\cK}(t)\te dE_R(z)\]
and
\[f(S)=\int_\Gamma \int_0^{2\pi}f_o(q^ls,t)dE_{\cK}(t)\te dE_S(z)\ ,\]
where the  function $f_o$ is complex-valued. Thus Theorem 
 \ref{hip2} reduces to the already proved "scalar" Theorem \ref{hip}.

We proceed to prove that really  
 for any $r,s\in\R_+$ and  $k,l=0,1,\dots,N-1$ we have
\beq
f(q^kr)f(q^ls)=f(q^ls)f(q^kr)\ 
\label{przem2}\ .
\eeq
Observe that operators $\faz R$ and  $\faz S$  have in a certain basis  
 the following matrix representation
\beq
\faz R=\left[\ba{ccccc}1&0&0&...&0\\0&q&0&...&0\\0&0&q^2&...&0\\...&...
&...&...&...\\
0&0&0&...&q^{N-1}\ea\right]
\label{frmac}
\eeq
and
\beq
\label{fsmac}
\faz S=\left[\ba{ccccc}0&1&0&...&0\\0&0&1&...&0\\0&0&0&...&1\\
...&...&...&...&...\\0&0&0&...&0  \ea\right]\ .
\eeq
 Let $T=S^{-1}R$. 
Then
\[
\faz 
T=q^{-\frac{1}{2}} (\faz S)^* \faz R \ ,
\]
so 
\beq
\label{ftmac}
\faz T=q^{-\frac{1}{2}}\left[\ba{ccccc}0&0&...&0&q^{N-1}\\1&0&...&0&0\\
0&q&...&0&1\\
...&...&...&...&...\\1&0&...&q^{N-2}&0  \ea\right]\ .
\eeq
We use the notion of generalized eigenvectors.
It is well known that a selfadjoint operator with continuous spectrum 
 acting on $\cH$ does not have eigenvectors. Still one can show that in the general case 
the generalized eigenvectors are continuous linear functionals 
on a certain dense locally convex subspace $\Phi\subset \cH$, provided with a much 
stronger topology than $\cH$. Then we get the same formulas as
 for discreet spectrum provided we replace scalar product by the duality relation  between $\Phi$ and $\Phi^\prime$.   
 This will be explained 
 by the example below, for general considerations see \cite{maueig}.

\bpr 
Let $\cH=\Lkr$
and
\beq
\label{np2}
 |R|=e^{\po} \mi |S|=e^{\qo} \mi |T|=\eh |S|^{-1}|R|=e^{\po-\qo}\ .\eeq
These operators have continuous spectra, so they 
do not have eigenvectors.
 There are however tempered distributions on $\R$ 
such that for every function $f$
  from the  Schwartz space of smooth functions on $\R$ 
 decreasing rapidly at infinity  $\sz$ we have 
\beq
\langle f|\;|R|\;|\Omega_r\rangle=r\langle f|\Omega_r\rangle\mi 
\langle f|\;|S|\;|\Phi_s\rangle=s\langle f|\Phi_s \rangle  \mi 
\langle f |\;|T|\;|\Psi_t\rangle =t\langle f|\Psi_t\rangle . 
\label{eigenvectors}\eeq
Such  $|\Omega_r\rangle$, $|\Phi_s\rangle$ and $|\Psi_t\rangle$ 
 are called generalized eigenvectors of 
 operators  $|R|$,$|S|$ and $|T|$ with  eigenvalues 
respectively $r$, $s$ and $t$.

An example of  generalized eigenvectors of
 operators \rf{np2} is
\[|\Omega_r\rangle=\frac{1}{\sqrt{2\pi\hbar}}e^{\frac{i}{\hbar}x\log r}\mi
|\Phi_s\rangle=\delta(\log s-x)\mi
|\Psi_t\rangle=\frac{1}{\sqrt{2\pi\hbar}}e^{i\frac{x^2}{2\hbar}}
e^{i\frac{x\log t}{\hbar}}\ .\]

Moreover, we will use notation of a type 
$\langle \Omega_r|\Phi_s\rangle$.

It should be understood in the following way: for any $f\in \sz$ we have
\beq\langle \Omega_r|f\rangle=\int_{\R} \langle \Omega_r|\Phi_s\rangle 
\langle \Phi_s|f\rangle ds\ .
\label{dlugie}\eeq
To shorten notation from now on we skip  
 the integration symbol, i.e. we write 
\[\langle \Omega_r|f\rangle=\langle \Omega_r|\Phi_s\rangle 
\langle \Phi_s|f\rangle  .\]
instead of \rf{dlugie}
The generalized eigenvectors $\Omega$ are said to have     
 the  {\em Dirac $\delta$ normalization} if  
\[\langle \Omega_r|\Omega_s\rangle =\delta (r-s)\ ,\]
where $\delta$ is the Dirac  $\delta$ distribution.
Note that generalized eigenvectors $\Omega$, $\Phi$ and $\Psi$ given above 
have  the   Dirac $\delta$ normalization.
\epr

Let $|\Omega_r\rangle$ be a generalized eigenvector of $R$  
with real eigenvalue $r$ and with  Dirac delta normalization.
Analogously, let  $|\Phi_s\rangle$ and $|\Psi_t\rangle$ denote 
generalized eigenvectors of 
 operators  $S$ and $T$ with real eigenvalues respectively  $s$ and $t$
 and with the Dirac delta normalization.

Let us define  $|\Omega_{k,r}\rangle$ using   
  the  vector $|\Omega_{r}\rangle$
 introduced above
\[|\Omega_{k,r}\rangle=e_k\te |\Omega_{r}\rangle\]
where and $k=0,1,\dots,N-1$ and 
\[e_0=\left[\ba{c}1\\0\\0\\ \vdots\\0\ea\right],
\;\;
e_1=\left[\ba{c}0\\1\\0\\\vdots\\0\ea\right],\;\;\dots,\;\;
e_{N-1}=\left[\ba{c}0\\0\\0\\\vdots\\1\ea\right].\ \]
Similarly, we define the vector  $|\Phi_{l,s}\rangle$ 
using the vector  $|\Phi_{s}\rangle$ introduced above
\[|\Phi_{l,s}\rangle=f_l\te |\Phi  _{s}\rangle\]
where \[f_l=\frac{1}{\sqrt{N}}\left[\ba{c}1\\q^l\\q^{2l}\\\vdots\\
q^{(N-1)l}\ea\right]\]
 and  $l=0,1,\dots,N-1$.

Analogously, the vector  $|\Psi_{m,t}\rangle$ 
 is given by 
\[|\Psi_{m,t}\rangle=g_m\te |\Psi  _{t}\rangle\]
where \[
(g_m)_p=\frac{1}{\sqrt{N}}q^{\frac{1}{2}(p^2-2p(m+1))}
\]
and  $m,
p=0,1,\dots,N-1$.

Note that $|\Omega_{k,r}\rangle$, $|\Phi_{l,s}\rangle$ and 
$|\Psi_{m,t}\rangle$ 
 are generalized eigenvectors of  $R$, $S$ and $T$, respectively,
 corresponding to generalized eigenvalues $q^kr$,$q^ls$ and $ q^mt$.

In order to prove  \rf{przem2} we compute matrix elements   
\[\langle \Omega_{k,r}|f(R)f(S)|\Phi_{l,s}\rangle\ .\]
Since the function $f$ satisfies exponential equation,  is Borel 
 and by  \rf{r+s}
 \[f(R)f(S)=f(R\dot{+}S)=F_N(T)^* f(S)F_N(T)\ .\]
Therefore
\[\langle \Omega_{k,r}|f(R)f(S)|\Phi_{l,s}\rangle =
\langle\Omega_{k,r}|F_N(T)^*|\Psi_{m,t}\rangle
\langle\Psi_{m,t}|f(S)|\Phi_{n,\tilde{s}}\rangle
\langle \Phi_{n,\tilde{s}}|F_N(T)|\Psi_{o,\tilde{t}}\rangle
\langle\Psi_{o,\tilde{t}}|\Phi_{l,s}\rangle
\ .\]
Hence
\[f(q^kr)f(q^ls)=\langle \Omega_{k,r}|\Phi_{l,s}\rangle^{-1}\langle\Omega_{k,r}|F_N(T)^*|\Psi_{m,t}\rangle
\langle\Psi_{m,t}|f(S)|\Phi_{n,\tilde{s}}\rangle
\langle \Phi_{n,\tilde{s}}|F_N(T)|\Psi_{o,\tilde{t}}\rangle
\langle\Psi_{o,\tilde{t}}|\Phi_{l,s}\rangle
\ .\]
It is easily checked that 
\[\langle \Omega_{k,r}|\Phi_{l,s}\rangle=
q^{kl}\langle\Omega_{r}|\Phi_{s}\rangle\]
\[\langle 
\Omega_{k,r}|\Psi_{m,t}\rangle=\frac{1}{\sqrt{N}}
q^{\frac{1}{2}(k^2-2k(m+1))}
\langle 
\Omega_{r}|\Psi_{t}\rangle
\]
\[\langle\Psi_{m,t}|\Phi_{n,\tilde{s}}\rangle=
\langle\Psi_{t}|\Phi_{\tilde{s}}\rangle
\frac{1}{N}\left(\sum_{p=o}^{N-1}
q^{p(m+n+1-\frac{p}{2})}\right)\]
As we see from the last formula, one has to compute 
 for any integer  $\alpha$  a sum
\[\sum_{p=o}^{N-1}
e^{\frac{2\pi i}{N}p(\alpha-\frac{p}{2})}=
e^{\frac{\pi i}{N}\alpha^2}
\sum_{p=o}^{N-1}
e^{-\frac{\pi i}{N}(p-\alpha)^2}=e^{\frac{\pi i}{N}\alpha^2}
\sum_{p=o}^{N-1}
e^{-\frac{\pi i}{N}p^2}\ .\]
In the appendix we derive the formula
 \beq
\label{suma2}
\sum_{p=o}^{N-1}
e^{-\frac{\pi i}{N}p^2}=\sqrt{N}e^{-\frac{\pi i}{4}}
\eeq
Hence
\[\sum_{p=o}^{N-1}
e^{\frac{2\pi i}{N}p(\alpha-\frac{p}{2})}=\sqrt{N}e^{\frac{\pi i}{N}\alpha^2}
e^{-\frac{\pi i}{4}}\]
so 
\[\langle\Psi_{m,t}|\Phi_{n,\tilde{s}}\rangle=
\frac{1}{\sqrt{N}}e^{-\frac{\pi i}{4}}
e^{\frac{\pi i}{N}(m+n+1)^2}
\langle\Psi_{t}|\Phi_{\tilde{s}}\rangle\]
Using Proposition 1.1  \cite{az+b}
  one can easily prove that
\bst
\label{wlFN}
For any  $m=0,1,\dots,N-1$ and $t \in R_+$ we have
\beq
\label{FNbar}
\ov{F_N(q^mt)}=e^{\frac{\pi i}{6}(\frac{2}{N}+\frac{N}{2})}
e^{-\frac{\pi i}{N}(m+1)^2}
e^{\frac{i}{2\hb}\log ^2t}F_N(q^{-m-2}t^{-1})\ ,
\eeq
\est
By Proposition above
\[\langle \Omega_{k,r}|F_N(T)^*|\Phi_{n,\tilde{s}}\rangle=
e^{\frac{\pi i}{6}(\frac{2}{N}+\frac{N}{2})}
e^{-\frac{\pi i}{N}(m+1)^2}
e^{\frac{i}{2\hb}\log ^2t}F_N(q^{-m-2}t^{-1})
\langle \Omega_{k,r}|\Psi_{m,t}\rangle\langle\Psi_{m,t}
|\Phi_{n,\tilde{s}}\rangle
=\]
\[=
e^{\frac{\pi i}{6}(\frac{2}{N}+\frac{N}{2})}
e^{-\frac{\pi i}{N}(m+1)^2}
e^{\frac{i}{2\hb}\log ^2t}F_N(q^{-m-2}t^{-1})
\frac{e^{-\frac{\pi i}{4}}}{N}
e^{\frac{i\hb}{2}(k^2-2k -2km +m^2+n^2+2mn +2m+2n+1)}
\langle 
\Omega_{r}|\Psi_{t}\rangle\langle\Psi_{t}|\Phi_{\tilde{s}}\rangle=\]
\[=\frac{e^{-\frac{\pi i}{4}}}{N}e^{\frac{\pi i}{6}(\frac{2}{N}+\frac{N}{2})}
e^{\frac{i}{2\hb}\log ^2t}
e^{\frac{i\hb}{2}(k^2-2k -2km +n^2+2mn +2n)}F_N(q^{-m-2}t^{-1})
\langle 
\Omega_{r}|\Psi_{t}\rangle\langle\Psi_{t}|\Phi_{\tilde{s}}\rangle\]
Inserting  $-m-2$ in the place of $m$, we get
\[\langle \Omega_{k,r}|F_N(T)^*|\Phi_{n,\tilde{s}}\rangle=\frac{e^{-\frac{\pi i}{4}}}{N}e^{\frac{\pi i}{6}(\frac{2}{N}+\frac{N}{2})}
e^{\frac{i}{2\hb}\log ^2t}
e^{\frac{i\hb}{2}(k^2+2k +2km +n^2-2mn -2n)}F_N(q^{m}t)
\langle 
\Omega_{r}|\Phi_{\tilde{s}}\rangle \ .\]

Observed that  by formula 1.36  \cite{qexp}, for any  $t\in\R$
\beq
\label{vts}
\overline{\vt(\log t)}=e^{-i\frac{\pi}{4}}c_\hbar^\prime e^{i\frac{\log^2 t}{2\hb}}
\vt(-\log t)
\ ,\eeq
where $ c_\hbar^\prime=
e^{i(\frac{\pi}{4}+\frac{\hb}{24}+\frac{\pi^2}{6\hb})}$. 

Moreover, we have prove in \cite{paper1} that 
\beq
\label{ef}
\langle\Omega_r|\vt(\log T)^*|\Phi_{s}\rangle=c_\hbar^\prime  
e^{-i\frac{\log^2 s}{\hbar}}\langle\Phi_{s}|\vt(\log T)|\Phi_r\rangle
\ .
\eeq
By \rf{vts} and \rf{ef}
\[e^{-\frac{\pi i}{4}}e^{\frac{i}{2\hb}\log ^2t}
F_N(q^mt^{-1})\langle 
\Omega_{r}|\Psi_{t}\rangle\langle\Psi_{t}|\Phi_{\tilde{s}}\rangle=
e^{-i\frac{\log^2\tilde{s}}{\hb}}F_N(q^mt)
\langle\Phi_{\tilde{s}}|\Psi_{t}\rangle\langle\Psi_{t}|\Phi_r\rangle\]
Moreover
\[ 
\frac{1}{N} F_N(q^m t)e^{\frac{\pi i}{N}(k^2+2km-2mn +2k +n^2 -2n)}=
e^{i\hb n^2}
\langle f_n|g_m\rangle F_N(q^m t)\langle g_m|f_k \rangle\]
Hence
\beq
\label{FNbarmac}
\langle \Omega_{k,r}|F_N(T)^*|\Phi_{n,\tilde{s}}\rangle=
e^{\frac{\pi i}{6}(\frac{2}{N}+\frac{N}{2})}
e^{i\hb n^2} e^{-i\frac{\log^2\tilde{s}}{\hb}}
\langle 
\Phi_{n,\tilde{s}}|F_N(T)|\Phi_{k,r}\rangle\  .
\eeq

Therefore
\[f(q^kr)f(q^ls)=q^{-kl}f(q^n\tilde{s})\langle\Omega_{r}|\Phi_{s}\rangle^{-1}
\langle \Omega_{k,r}|F_N(T)^*|\Phi_{n,\tilde{s}}\rangle
\langle\Phi_{n,\tilde{s}}|F_N(T)|\Phi_{l,s}\rangle=\]
\[=e^{\frac{\pi i}{6}(\frac{2}{N}+\frac{N}{2})}e^{i\hb n^2}
 e^{-i\frac{\log^2\tilde{s}}{\hb}}
f(q^n\tilde{s})
q^{-kl}\langle\Omega_{r}|\Phi_{s}\rangle^{-1}
\langle\Phi_{n,\tilde{s}}|F_N(T)|\Phi_{k,r}\rangle
\langle\Phi_{n,\tilde{s}}|F_N(T)|\Phi_{l,s}\rangle\]

The expression below is clearly  symmetric with respect to 
 swapping  $k\leftrightarrow l$ with $r\leftrightarrow s$, 
which completes the proof.\hfill\qed

\appendix
\section{Deriving  formula  \rf{suma2}}
\label{dodsuma}
We will derive the formula
\beq
\label{suma}
\sum_{p=0}^{N-1}e^{\frac{\pi i}{N}p^2},
\eeq
where  $N$ is a non-zero  even number.

Let us introduce the notation
\[a_p=e^{\frac{\pi i}{N}p^2}\mi S_N=\sum_{p=0}^{N-1}a_p\ .\]

\begin{figure}
\centerline{\hfill
\epsfxsize=100mm
\epsfysize=100mm
\epsffile{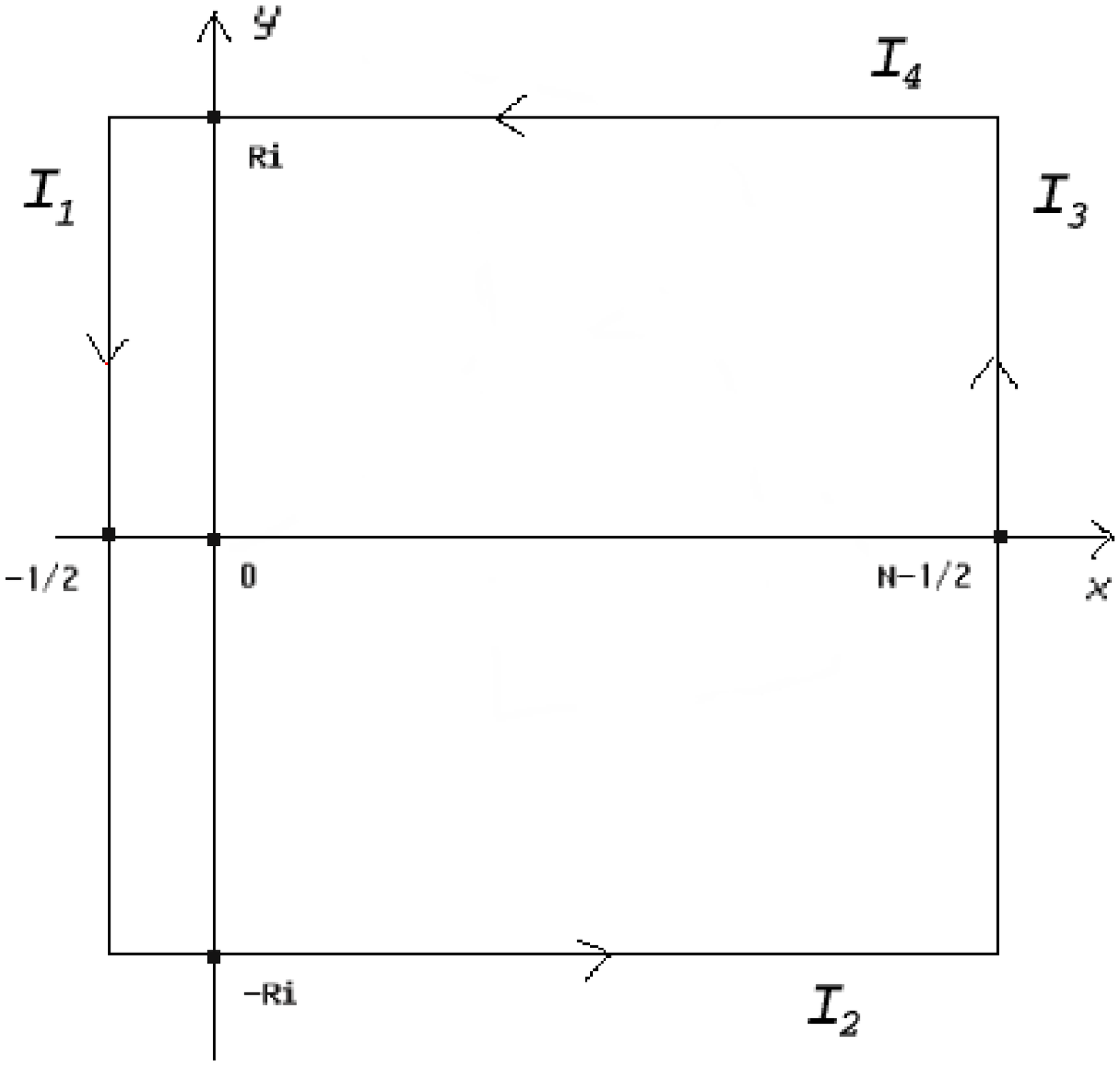}
\hfill
}
\caption{\label{kontur} contour $\Gamma$}
\end{figure}

In order to compute \rf{suma} we integrate the function  
\beq
\label{function}
f(z)=\frac{e^{\frac{\pi i}{N}z^2}}
{e^{2\pi iz}-1}
\eeq
over the contour $\Gamma$ as on the Figure \ref{kontur}.

Since $N$ is an even number, it follows that for any  $z\in \C$ 
 from the domain of  $f$ 
we have
\beq
\label{i1+i3}
f(z+N)-f(z)=e^{\frac{\pi i z^2}{N}}\ .
\eeq
Let us introduce the notation
\[I_1=\int_{R}^{-R}
f(-\frac{1}{2}+iy)idy
\mi
I_2=\int_{-\frac{1}{2}}^{N-\frac{1}{2}}f(x-iR)dx\]
\[I_3=\int_{-R}^{R}
f(N-\frac{1}{2}+iy)
idy
\mi
I_4=\int_{N-\frac{1}{2}}^{-\frac{1}{2}}
f(x+iR)dx
\]
Observe that
\beq
\label{wres1}
\int_\Gamma f(z)dz=I_1+I_2+I_3+I_4
\eeq
Function  $f$ in the interior of  $\Gamma$  has simple poles at the  
 points $z=1,2,\dots,\frac{N}{2}$, and  residues at these points 
 are
\[\res_{z=p} f(z)=\frac{a_p}{2\pi i}\ .\] 
Therefore by the Residue Theorem
\beq
\label{wres2}
\int_\Gamma f(z)dz=2\pi i\sum_{p=0}^{N-1}\res_{z=p} f(z)=
a_0+a_2+\dots +a_{N-1}=S_{N}\ .
\eeq
Comparing \rf{wres1} and  \rf{wres2} we get
\beq
\label{wres4}
S_N=I_1+I_2+I_3+I_4\ .
\eeq
We proceed now to calculate integrals  $I_2$ and $I_4$.
\[I_2=\int_{-\frac{1}{2}}^{N-\frac{1}{2}}\frac{e^{\frac{\pi i}{N}(x-iR)^2}}
{e^{2\pi i(x-iR)}-1}dx=\int_{-\frac{1}{2}}^{N-\frac{1}{2}}
\frac{e^{\frac{\pi i}{N}(x^2-R ^2)}e^{\frac{2\pi }{N}R x}}
{e^{2\pi (R+ix)}-1}dx\]
Since
\[ |e^{\frac{\pi i}{N}(x^2-R ^2)}e^{\frac{2\pi }{N}R x}|
\leq 
e^{\frac{2\pi }{N}R x}
\mi |e^{2\pi (R+ix)}-1|\geq e^{2\pi R}-1\ ,\]
it follows that 
\[I_2<\left|\int_{-\frac{1}{2}}^{N-\frac{1}{2}}
\frac{e^{\frac{\pi i}{N}(x-iR)^2}}
{e^{2\pi (R+ix)}-1}dx\right|\leq
\frac{1}
{e^{2\pi R}-1}
\int_{-\frac{1}{2}}^{N-\frac{1}{2}}
e^{\frac{2\pi }{N}Rx}
dx=\frac{Ne^{-\pi R}}{2\pi R}
\frac{e^{2\pi R}-1}
{e^{2\pi R}-1}
=\frac{Ne^{-\pi R}}{2\pi R}\ .\]
Hence $I_2$  converges to 0, when  $R$ goes to $+\infty$.
Using similar estimation one can prove that   
$I_4$ converges also  to 0, when $R$  goes to  $+\infty$.
Moreover, using formula  \rf{i1+i3} we get
\[I_1+I_3 =2\int_{-R}^{R}\left(f(N-\frac{1}{2}+it)-
f(-\frac{1}{2}+it)\right)
idy=2\int_{-R}^{R}e^{\frac{\pi i (it-\frac{1}{2})^2}{N}}
idy
=2i\int_{-R}^{R}e^{\frac{-\pi i}{N}(t+\frac{1}{2}i)^2 }
dy\ .\]
When  $R\rightarrow +\infty$, then $ I_1+I_3$ takes form 
\[2i\int_{-\infty}^{+\infty}e^{\frac{-\pi i}{N}(t+\frac{1}{2}i)^2 }
dy=2i\frac{\sqrt{N}}{\sqrt{\pi}}\int_{-\infty}^{+\infty}e^{ -i (t+\frac{1}{2}i)^2 }
dy\]
After passing to the limit  $R\rightarrow +\infty$ into  
 the formula \rf{wres4} we get  
\[S_N=2i\frac{\sqrt{N}}{\sqrt{\pi}}\int_{-\infty}^{+\infty}
e^{ -i (t+\frac{1}{2}i)^2 }
dy\]
It is easily seen that
\[S_2=1+i\ .\]
Hence
\[\int_{-\infty}^{+\infty}e^{ -i y^2 }
dy=\frac{\sqrt{\pi}}{2\sqrt{2}i}(i+1)=
-i\frac{\sqrt{\pi}}{2}e^{\frac{\pi i}{4}}\ .\]
Thus for any even $N\geq 2$  we have 
\[\sum_{p=0}^{N-1}e^{\frac{\pi i}{N}p^2}=\sqrt{N}e^{\frac{\pi i}{4}}\]
\section*{Acknowledgments}
This is part of the author's Ph.D. thesis \cite{phd}, written under the supervision of Professor Stanis{\l}aw L. Woronowicz at the Department of Mathematical Methods in 
Physics, Warsaw University. The author is greatly  
indebted to Professor Stanis{\l}aw L. Woronowicz for stimulating discussions and important hints 
 and comments. The author also wishes to thank 
Professor Wies{\l}aw Pusz and Professor Marek  Bo\.zejko for
several  helpful suggestions and Piotr So\l tan for 
reading carefully the manuscript.

\end{document}